# PRICING AND TRADING CREDIT DEFAULT SWAPS IN A HAZARD PROCESS MODEL

By Tomasz R. Bielecki,[1] Monique Jeanblanc[2]
and Marek Rutkowski[3]

*Illinois Institute of Technology, Université d'Évry Val d'Essonne and
University of New South Wales and Warsaw University of Technology*

In the paper we study dynamics of the arbitrage prices of credit default swaps within a hazard process model of credit risk. We derive these dynamics without postulating that the immersion property is satisfied between some relevant filtrations. These results are then applied so to study the problem of replication of general defaultable claims, including some basket claims, by means of dynamic trading of credit default swaps.

**1. Introduction.** An inspection of the existing literature in the area of credit risk shows that the vast majority of papers focus on the risk-neutral valuation of credit derivatives without even mentioning the issue of hedging. This is somewhat surprising since, as is well known, the major argument supporting the risk-neutral valuation is the existence of hedging strategies for attainable contingent claims. In this paper we shall deal with the credit default swaps market only. Valuation formulae for credit derivatives traded on the CDS market are provided, for instance, in Brigo [11], Brigo and Morini [12], Hull and White [17], Schönbucher [24] and Wu [26], who deal with different products and/or models. There also exists a slowly growing number of papers in which the issue of hedging of defaultable claims is analyzed in a more systematic way; to mention a few: Arvanitis and Laurent [2], Bélanger, Shreve and Wong [3], Bielecki, Jeanbanc and Rutkowski [5, 6,

Received August 2007; revised January 2008.

[1]Supported in part by NSF Grants 0202851 and 0604789, and by Moody's Corporation Grant 5-55411.

[2]Supported in part by Itô33, and Moody's Corporation Grant 5-55411.

[3]Supported in part by ARC Discovery Project Grant DP0881460.

7, 8], Blanchet–Scalliet and Jeanblanc [10], Collin–Dufresne and Hugonnier [13], Frey and Backhaus [16], Kurtz and Riboulet [20], Laurent [21] and Laurent et al. [22]. In particular, in [5, 7] we studied hedging of defaultable claims in a general, semi-martingale setting. On the other hand, the pilot paper [8] studies problems analogous to those studied here, but under the assumption that the reference filtration (denoted as $\mathbb{F}$ in what follows) is a trivial filtration.

From the practical perspective, it is common to split the risk of a credit derivative into three components: the *default risk* (that is, the jump risk associated with some particular credit event), the *spread risk* (i.e., the risk due to the volatile character of the pre-default values of a credit) and the *correlation risk* due to interdependence of the underlying credit names. The pertinent issue is thus to find a mathematically rigorous way of dealing simultaneously with all three kinds of credit risks. Our main results, Theorems 2.1 and 3.1, show that in a generic hazard process model driven by a Brownian motion it is, in principle, possible to perfectly hedge all sorts of risks in a unified manner, provided that a large enough number of liquid CDSs are traded. Specifically, one set of conditions address the issue of, hedging default risk, while other ones allow us to effectively deal with spread and correlation risks. Of course, all these conditions need to be simultaneously satisfied for a perfect hedging. It is clear from these formulae that hedging of the default risk relies on keeping under control the unexpected jumps that may come as a surprise at any moment and that are modeled by pure jump martingales. By contrast, hedging of spread and correlation risks hinges on more standard techniques related to volatilities and correlations of underlying continuous driving martingales. Let us finally note that in our previous paper Bielecki, Jeanblanc and Rutkowski [7] we have shown that in the case of a *survival claim* (i.e., a defaultable claim with zero recovery a default), it is enough to focus on hedging of the spread risk, provided that hedging instruments are also subject to the zero recovery scheme. We shall work throughout within the so-called hazard process (or reduced-form) approach, as opposed to the structural approach in which hedging can be dealt with using the classic Black–Scholes-like approach. For general results within this methodology, we refer the reader to, among others, Bélanger, Shreve and Wong [3], Bielecki and Rutkowski [4], Elliott, Jeanblanc and Yor [15], Jeanblanc and Le Cam [18], Jeanblanc and Rutkowski [19] and Schönbucher and Schubert [25].

Our program can be summarized as follows. We start by deriving the risk-neutral dynamics for the prices of defaultable claims. Next, we show that the risk-neutral pricing of defaultable claims (such as, credit default options and first-to-default swaps) can be supported through replication of these claims by dynamic trading of a suitable family of single-name credit default swaps. In Section 2 we address the issue of valuation and hedging of defaultable claims in the market with traded CDSs with different maturities but with the same reference credit name. Results of this section may thus be applied, for instance, to a single-name credit default swaption. We show that replication of defaultable claims can be done using either a family of CDSs with fixed spreads and different maturities or using the associated family of virtual market CDSs, which may serve as the proxy for the market CDSs (in practice, CDSs are issued on a daily basis at the current market spread). In Section 3 we first derive the dynamics for a family of single-name CDSs in the case of several correlated credit names. Obviously, the fact that the default on a particular name occurs has an impact on the dynamics of CDS written on nondefaulted names. Subsequently, we extend some results obtained previously by Bielecki, Jeanblanc and Rutkowski [8] in the case of a trivial reference filtration to a more practically appealing case of a market model in which hazard rates are driven by a multidimensional Brownian motion. Let us admit that the results of this paper are of a rather abstract nature, in the sense that no explicit examples of hedging strategies for standard credit derivatives are presented; they will be studied in the follow-up work by Bielecki, Jeanblanc and Rutkowski [9].

**2. Single-name credit default swap market.** A strictly positive random variable $\tau$, defined on a probability space $(\Omega, \mathcal{G}, \mathbb{Q})$, is termed a *random time*. In view of its financial interpretation, we will refer to it as a *default time*. We define the default indicator process $H_t = \mathbb{1}_{\{\tau \leq t\}}$ and we denote by $\mathbb{H}$ the filtration generated by this process. We assume that we are given, in addition, some auxiliary filtration $\mathbb{F}$ and we write $\mathbb{G} = \mathbb{H} \vee \mathbb{F}$, meaning that we have $\mathcal{G}_t = \sigma(\mathcal{H}_t, \mathcal{F}_t)$ for every $t \in \mathbb{R}_+$. The filtration $\mathbb{G}$ is referred to as the *full filtration*. It is clear that $\tau$ is a $\mathbb{H}$-stopping time, as well as a $\mathbb{G}$-stopping time (but not necessarily an $\mathbb{F}$-stopping time). All processes are defined on the space $(\Omega, \mathbb{G}, \mathbb{Q})$, where $\mathbb{Q}$ is to be interpreted as the real-life (i.e., statistical) probability measure. Unless otherwise stated, all processes considered in what follows are assumed to be $\mathbb{G}$-adapted and with càdlàg sample paths.

2.1. *Price dynamics in a single-name model.* We assume that the underlying market model is arbitrage-free, meaning that it admits a *spot martingale measure* $\mathbb{Q}^*$ (not necessarily unique) equivalent to $\mathbb{Q}$. A *spot martingale measure* is associated with the choice of the savings account $B$ as a numéraire, in the sense that the price process of any tradeable security, which pays no coupons or dividends, is a $\mathbb{G}$-martingale under $\mathbb{Q}^*$, when it is discounted by the *savings account* $B$. As usual, $B$ is given by

$$(1) \qquad B_t = \exp\left(\int_0^t r_u \, du\right) \qquad \forall t \in \mathbb{R}_+,$$

where the short-term $r$ is assumed to follow an $\mathbb{F}$-progressively measurable stochastic process. The choice of a suitable term structure model is arbitrary and it is not discussed in the present work.

Let us denote by $G_t = \mathbb{Q}^*(\tau > t \mid \mathcal{F}_t)$ the *survival process* of $\tau$ with respect to a filtration $\mathbb{F}$. We postulate that $G_0 = 1$ and $G_t > 0$ for every $t \in \mathbb{R}_+$ (hence, the case where $\tau$ is an $\mathbb{F}$-stopping time is excluded) so that the *hazard process* $\Gamma = -\ln G$ of $\tau$ with respect to the filtration $\mathbb{F}$ is well defined.

For any $\mathbb{Q}^*$-integrable and $\mathcal{F}_T$-measurable random variable $Y$, the following classic formula holds (see, e.g., Chapter 5 in [4] or [19]):

$$\tag{2} \mathbb{E}_{\mathbb{Q}^*}(\mathbb{1}_{\{T < \tau\}} Y \mid \mathcal{G}_t) = \mathbb{1}_{\{t < \tau\}} G_t^{-1} \mathbb{E}_{\mathbb{Q}^*}(G_T Y \mid \mathcal{F}_t).$$

Clearly, the process $G$ is a bounded $\mathbb{G}$-supermartingale and thus it admits the unique Doob–Meyer decomposition $G = \mu - \nu$, where $\mu$ is a martingale part and $\nu$ is a predictable increasing process. Note that if $G$ is continuous, then the processes $\mu$ and $\nu$ are continuous as well.

We shall work throughout under the following standing assumption.

ASSUMPTION 2.1.   We postulate that $G$ is a continuous process and the increasing process $\nu$ in its Doob–Meyer decomposition is absolutely continuous with respect to the Lebesgue measure, so that $d\nu_t = \upsilon_t \, dt$ for some $\mathbb{F}$-progressively measurable, nonnegative process $\upsilon$. We denote by $\lambda$ the $\mathbb{F}$-progressively measurable process defined as $\lambda_t = G_t^{-1} \upsilon_t$.

Let us note for the further reference that under Assumption 2.1 we have $dG_t = d\mu_t - \lambda_t G_t \, dt$, where the $\mathbb{F}$-martingale $\mu$ is continuous. Moreover, in view of the Lebesgue dominated convergence theorem, continuity of $G$ implies that the expected value $\mathbb{E}_{\mathbb{Q}^*}(G_t) = \mathbb{Q}^*(\tau > t)$ is a continuous function, and thus, $\mathbb{Q}^*(\tau = t) = 0$, for any fixed $t \in \mathbb{R}_+$. Finally, it is known (see, e.g., Lemma 3.2 in [15], or [19]) that under Assumption 2.1 the process $M$, given by

$$\tag{3} M_t = H_t - \Lambda_{t \wedge \tau} = H_t - \int_0^{t \wedge \tau} \lambda_u \, du = H_t - \int_0^t (1 - H_u) \lambda_u \, du,$$

is a $\mathbb{G}$-martingale, where the increasing, absolutely continuous, $\mathbb{F}$-adapted process $\Lambda$ is given by

$$\tag{4} \Lambda_t = \int_0^t G_u^{-1} \, d\nu_u = \int_0^t \lambda_u \, du.$$

The $\mathbb{F}$-progressively measurable process $\lambda$ is called the *default intensity* with respect to $\mathbb{F}$.

REMARK.   Results of this paper can be extended to the case where $G$ is not assumed to be continuous, under the assumptions that $\tau$ avoids the $\mathbb{F}$-stopping times and that $\nu$ is a continuous process. The continuity of $\nu$ is required to obtain a continuous compensator of $H$, that is, to work with a totally inaccessible $\mathbb{G}$-stopping time $\tau$. Note also that under the assumption that $\nu$ is absolutely continuous with respect to the Lebesgue measure, the hypothesis on continuity of $G$ is not needed. Indeed, given that $d\nu_t = \upsilon_t \, dt$, the compensator of $H$ then equals $\Lambda_{t\wedge\tau} = \int_0^{t\wedge\tau}(G_{u-})^{-1}\upsilon_u \, du$, and thus we also have that

$$\Lambda_{t\wedge\tau} = \int_0^{t\wedge\tau}(G_u)^{-1}\upsilon_u \, du = \int_0^{t\wedge\tau}\lambda_u \, du.$$

2.1.1. *Defaultable claims.* We are in the position to introduce the concept of a defaultable claim. Of course, we work here within a single-name framework, so that $\tau$ is the moment of default of the reference credit name.

DEFINITION 2.1. By a *defaultable claim maturing at $T$*, we mean the quadruple $(X, A, Z, \tau)$, where $X$ is an $\mathcal{F}_T$-measurable random variable, $A = (A_t)_{t\in[0,T]}$ is an $\mathbb{F}$-adapted, continuous process of finite variation with $A_0 = 0$, $Z = (Z_t)_{t\in[0,T]}$ is an $\mathbb{F}$-predictable process, and $\tau$ is a random time.

The financial interpretation of components of a defaultable claim becomes clear from the following definition of the *dividend process $D$*, which describes all cash flows associated with a defaultable claim over its lifespan $]0, T]$, that is, after the contract was initiated at time 0 (of course, the choice of 0 as the inception date is merely a convention). The dividend process might have been called the *total cash flow process*; we have chosen the term "dividend process" for the sake of brevity.

DEFINITION 2.2. The *dividend process $D = (D_t)_{t\in\mathbb{R}_+}$* of the above defaultable claim maturing at $T$ equals, for every $t \in \mathbb{R}_+$,

$$D_t = X\mathbb{1}_{\{T<\tau\}}\mathbb{1}_{[T,\infty[}(t) + \int_{]0,t\wedge T]}(1-H_u)\,dA_u + \int_{]0,t\wedge T]}Z_u\,dH_u.$$

It is clear that the dividend process $D$ is a process of finite variation on $[0,T]$. The financial interpretation of $D$ is as follows: $X$ is the *promised payoff*, $A$ represents the process of *promised dividends* and the process $Z$, termed the *recovery process*, specifies the recovery payoff at default. It is worth stressing that, according to our convention, the cash payment (premium) at time 0 is not included in the dividend process $D$ associated with a defaultable claim.

2.1.2. *Price dynamics of a defaultable claim.* For any fixed $t \in [0, T]$, the process $D_u - D_t, u \in [t, T]$, represents all cash flows from a defaultable claim received by an investor who purchased it at time $t$. Of course, the process $D_u - D_t$ may depend on the past behavior of the claim as well as on the history of the market prior to $t$. The past dividends are not valued by the market, however, so that the current market value at time $t \in [0, T]$ of a defaultable claim (i.e., the price at which it trades at time $t$) reflects only future cash flows to be paid/received over the time interval $]t, T]$. This leads to the following definition of the ex-dividend price of a defaultable claim.

DEFINITION 2.3.   The *ex-dividend price* process $S$ of a defaultable claim $(X, A, Z, \tau)$ equals, for every $t \in [0, T]$,

$$(5) \qquad S_t = B_t \mathbb{E}_{\mathbb{Q}^*} \left( \int_{]t,T]} B_u^{-1} \, dD_u \, | \, \mathcal{G}_t \right).$$

Obviously, $S_T = 0$ for any dividend process $D$. We work throughout under the natural integrability assumptions,

$$\mathbb{E}_{\mathbb{Q}^*} |B_T^{-1} X| < \infty, \qquad \mathbb{E}_{\mathbb{Q}^*} \left| \int_{]0,T]} B_u^{-1} (1 - H_u) \, dA_u \right| < \infty,$$

$$\mathbb{E}_{\mathbb{Q}^*} |B_{\tau \wedge T}^{-1} Z_{\tau \wedge T}| < \infty,$$

which ensure that the ex-dividend price $S_t$ is well defined for any $t \in [0, T]$. We will later need the following technical assumption:

$$(6) \qquad \mathbb{E}_{\mathbb{Q}^*} \left( \int_0^T (B_u^{-1} Z_u)^2 \, d\langle \mu \rangle_u \right) < \infty.$$

We first derive a convenient representation for the ex-dividend price $S$ of a defaultable claim.

PROPOSITION 2.1.   *The ex-dividend price of the defaultable claim* $(X, A, Z, \tau)$ *equals, for* $t \in [0, T[$,

$$(7) \quad S_t = \mathbb{1}_{\{t < \tau\}} \frac{B_t}{G_t} \mathbb{E}_{\mathbb{Q}^*} \left( B_T^{-1} G_T X + \int_t^T B_u^{-1} G_u (Z_u \lambda_u \, du + dA_u) \, | \, \mathcal{F}_t \right).$$

PROOF.   For any $t \in [0, T[$, the ex-dividend price is given by the conditional expectation

$$S_t = B_t \mathbb{E}_{\mathbb{Q}^*} \left( B_T^{-1} X \mathbb{1}_{\{T < \tau\}} + \int_{t \wedge \tau}^{T \wedge \tau} B_u^{-1} \, dA_u + B_\tau^{-1} Z_\tau \mathbb{1}_{\{t < \tau \leq T\}} \, | \, \mathcal{G}_t \right).$$

Let us fix $t$ and let us introduce two auxiliary processes $Y = (Y_u)_{u \in [t,T]}$ and $R = (R_u)_{u \in [t,T]}$ by setting

$$Y_u = \int_t^u B_v^{-1} \, dA_v, \qquad R_u = B_u^{-1} Z_u + \int_t^u B_v^{-1} \, dA_v = B_u^{-1} Z_u + Y_u.$$

Then $S_t$ can be represented as follows:

$$S_t = B_t \mathbb{E}_{\mathbb{Q}^*}(B_T^{-1} X \mathbb{1}_{\{T < \tau\}} + \mathbb{1}_{\{T < \tau\}} Y_T + R_\tau \mathbb{1}_{\{t < \tau \leq T\}} \mid \mathcal{G}_t).$$

We use directly formula (2) in order to evaluate the conditional expectations

$$B_t \mathbb{E}_{\mathbb{Q}^*}(\mathbb{1}_{\{T < \tau\}} B_T^{-1} X \mid \mathcal{G}_t) = \mathbb{1}_{\{t < \tau\}} \frac{B_t}{G_t} \mathbb{E}_{\mathbb{Q}^*}(B_T^{-1} G_T X \mid \mathcal{F}_t)$$

and

$$B_t \mathbb{E}_{\mathbb{Q}^*}(\mathbb{1}_{\{T < \tau\}} Y_T \mid \mathcal{G}_t) = \mathbb{1}_{\{t < \tau\}} \frac{B_t}{G_t} \mathbb{E}_{\mathbb{Q}^*}(G_T Y_T \mid \mathcal{F}_t).$$

In addition, we will use of the following formula (see, e.g., [4]):

$$(8) \qquad \mathbb{E}_{\mathbb{Q}^*}(\mathbb{1}_{\{t < \tau \leq T\}} R_\tau \mid \mathcal{G}_t) = -\mathbb{1}_{\{t < \tau\}} \frac{1}{G_t} \mathbb{E}_{\mathbb{Q}^*}\left( \int_t^T R_u \, dG_u \,\Big|\, \mathcal{F}_t \right),$$

which is known to be valid for any $\mathbb{F}$-predictable process $R$ such that $\mathbb{E}_{\mathbb{Q}^*}|R_\tau| < \infty$. We thus obtain, for any $t \in [0, T[$,

$$S_t = \mathbb{1}_{\{t < \tau\}} \frac{B_t}{G_t} \mathbb{E}_{\mathbb{Q}^*}\left( B_T^{-1} G_T X + G_T Y_T - \int_t^T (B_u^{-1} Z_u + Y_u) \, dG_u \,\Big|\, \mathcal{F}_t \right),$$

Moreover, since $dG_t = d\mu_t - \lambda_t G_t \, dt$, where $\mu$ is an $\mathbb{F}$-martingale, we also obtain

$$\mathbb{1}_{\{t < \tau\}} \frac{B_t}{G_t} \mathbb{E}_{\mathbb{Q}^*}\left( - \int_t^T B_u^{-1} Z_u \, dG_u \,\Big|\, \mathcal{F}_t \right)$$

$$= \mathbb{1}_{\{t < \tau\}} \frac{B_t}{G_t} \mathbb{E}_{\mathbb{Q}^*}\left( \int_t^T B_u^{-1} G_u Z_u \lambda_u \, du \,\Big|\, \mathcal{F}_t \right),$$

where we have used (6). To complete the proof, it remains to observe that $G$ is a continuous semimartingale and $Y$ is a continuous process of finite variation with $Y_t = 0$, so that the Itô integration by parts formula yields

$$G_T Y_T - \int_t^T Y_u \, dG_u = \int_t^T G_u \, dY_u = \int_t^T B_u^{-1} G_u \, dA_u,$$

where the second equality follows from the definition of $Y$. We conclude that (7) holds for any $t \in [0, T[$, as required. $\quad\square$

Formula (7) implies that the ex-dividend price $S$ satisfies, for every $t \in [0, T]$,

$$(9) \qquad\qquad\qquad S_t = \mathbb{1}_{\{t < \tau\}} \widetilde{S}_t$$

for some $\mathbb{F}$-adapted process $\widetilde{S}$, which is termed the *ex-dividend pre-default price* of a defaultable claim. Note that $S$ may not be continuous at time $T$, in which case $S_{T-} \neq S_T = 0$.

DEFINITION 2.4. The *cumulative price* process $S^c$ associated with the dividend process $D$ is defined by setting, for every $t \in [0, T]$,

$$(10) \qquad S_t^c = B_t \mathbb{E}_{\mathbb{Q}^*} \left( \int_{]0,T]} B_u^{-1} \, dD_u \,\Big|\, \mathcal{G}_t \right) = S_t + B_t \int_{]0,t]} B_u^{-1} \, dD_u.$$

Note that the discounted cumulative price $B^{-1} S^c$ is a $\mathbb{G}$-martingale under $\mathbb{Q}^*$. It follows immediately from (7) and (10) that the following corollary to Proposition 2.1 is valid.

COROLLARY 2.1. *The cumulative price of the defaultable claim* $(X, A, Z, \tau)$ *equals, for* $t \in [0, T]$,

$$S_t^c = \mathbb{1}_{\{t < \tau\}} \frac{B_t}{G_t} \mathbb{E}_{\mathbb{Q}^*} \left( B_T^{-1} G_T X \mathbb{1}_{\{t < T\}} + \int_t^T B_u^{-1} G_u (Z_u \lambda_u \, du + dA_u) \,\Big|\, \mathcal{F}_t \right)$$
$$+ B_t \int_{]0,t]} B_u^{-1} \, dD_u.$$

The *pre-default cumulative price* is the unique $\mathbb{F}$-adapted process $\widetilde{S}^c$ that satisfies, for every $t \in [0, T]$,

$$(11) \qquad \mathbb{1}_{\{t < \tau\}} S_t^c = \mathbb{1}_{\{t < \tau\}} \widetilde{S}_t^c.$$

Our next goal is to derive the dynamics under $\mathbb{Q}^*$ for (pre-default) prices and of a defaultable claim in terms of some $\mathbb{G}$-martingales and $\mathbb{F}$-martingales. To simplify the presentation, we shall work from now on under the following standing assumptions.

ASSUMPTION 2.2. We assume that all $\mathbb{F}$-martingales are continuous processes.

The following auxiliary result is well known (see, e.g., Lemma 5.1.6 in [4]). Recall that $\mu$ is the $\mathbb{F}$-martingale appearing in the Doob–Meyer decomposition of $G$.

LEMMA 2.1. *Let* $n$ *be any* $\mathbb{F}$-*martingale. Then the process* $\widehat{n}$ *given by*

$$(12) \qquad \widehat{n}_t = n_{t \wedge \tau} - \int_0^{t \wedge \tau} G_u^{-1} \, d\langle n, \mu \rangle_u$$

*is a continuous* $\mathbb{G}$-*martingale.*

In particular, the process $\widehat{\mu}$ given by

$$(13) \qquad \widehat{\mu}_t = \mu_{t \wedge \tau} - \int_0^{t \wedge \tau} G_u^{-1} \, d\langle \mu, \mu \rangle_u$$

is a continuous $\mathbb{G}$-martingale.

In the next result we deal with the dynamics of the ex-dividend price process $S$. Recall that the $\mathbb{G}$-martingale $M$ is given by (3).

PROPOSITION 2.2. *The dynamics of the ex-dividend price $S$ on $[0, T]$ are*

$$
\begin{aligned}
(14) \quad dS_t = {} & -S_{t-}\, dM_t + (1 - H_t)((r_t S_t - \lambda_t Z_t)\, dt - dA_t) \\
& + (1 - H_t) G_t^{-1} (B_t\, dm_t - S_t\, d\mu_t) \\
& + (1 - H_t) G_t^{-2} (S_t\, d\langle \mu \rangle_t - B_t\, d\langle \mu, m \rangle_t),
\end{aligned}
$$

*where the continuous $\mathbb{F}$-martingale $m$ is given by the formula*

$$
(15) \quad m_t = \mathbb{E}_{\mathbb{Q}^*}\left( B_T^{-1} G_T X + \int_0^T B_u^{-1} G_u (Z_u \lambda_u\, du + dA_u) \,\Big|\, \mathcal{F}_t \right).
$$

PROOF. We shall first derive the dynamics of the pre-default ex-dividend price $\widetilde{S}$. In view of (7), the price $S$ can be represented as follows, for $t \in [0, T[$:

$$
S_t = \mathbb{1}_{\{t < \tau\}} \widetilde{S}_t = \mathbb{1}_{\{t < \tau\}} B_t G_t^{-1} U_t,
$$

where the auxiliary process $U$ equals

$$
U_t = m_t - \int_0^t B_u^{-1} G_u Z_u \lambda_u\, du - \int_0^t B_u^{-1} G_u\, dA_u,
$$

where, in turn, the continuous $\mathbb{F}$-martingale $m$ is given by (15). It is thus obvious that $\widetilde{S} = BG^{-1}U$ for $t \in [0, T[$ (of course, $\widetilde{S}_T = 0$). Since $G = \mu - \nu$, an application of Itô's formula leads to

$$
\begin{aligned}
d(G_t^{-1} U_t) = {} & G_t^{-1}\, dm_t - B_t^{-1} Z_t \lambda_t\, dt - B_t^{-1}\, dA_t \\
& + U_t (G_t^{-3}\, d\langle \mu \rangle_t - G_t^{-2}\, (d\mu_t - d\nu_t)) - G_t^{-2}\, d\langle \mu, m \rangle_t.
\end{aligned}
$$

Therefore, since under the present assumptions $d\nu_t = \lambda_t G_t\, dt$, using again Itô's formula, we obtain

$$
\begin{aligned}
(16) \quad d\widetilde{S}_t = {} & ((\lambda_t + r_t)\widetilde{S}_t - \lambda_t Z_t)\, dt - dA_t + G_t^{-1}(B_t\, dm_t - \widetilde{S}_t\, d\mu_t) \\
& + G_t^{-2}(\widetilde{S}_t\, d\langle \mu \rangle_t - B_t\, d\langle \mu, m \rangle_t).
\end{aligned}
$$

Note that, under the present assumptions, the pre-default ex-dividend price $\widetilde{S}$ follows on $[0, T[$ a continuous process with dynamics given by (16). This means that $S_{t-} = \widetilde{S}_t$ on $\{t \le \tau\}$ for any $t \in [0, T[$. Moreover, since $G$ is continuous, we have that $\mathbb{Q}^*(\tau = T) = 0$. Hence, for the process $S_t = (1 - H_t)\widetilde{S}_t$ we obtain, for every $t \in [0, T]$,

$$
\begin{aligned}
(17) \quad dS_t = {} & -S_{t-}\, dM_t + (1 - H_t)((r_t S_t - \lambda_t Z_t)\, dt - dA_t) \\
& + (1 - H_t) G_t^{-1}(B_t\, dm_t - S_t\, d\mu_t) \\
& + (1 - H_t) G_t^{-2}(S_t\, d\langle \mu \rangle_t - B_t\, d\langle \mu, m \rangle_t).
\end{aligned}
$$

This finishes the proof of the proposition.  □

Let us now examine the dynamics of the cumulative price. As expected, the discounted cumulative price $B^{-1}S^c$ is a $\mathbb{G}$-martingale under $\mathbb{Q}^*$ [see formula (19) below].

COROLLARY 2.2.   *The dynamics of the cumulative price $S^c$ on $[0, T]$ are*

(18)
$$\begin{aligned} dS_t^c = {}& r_t S_t^c \, dt + (Z_t - S_{t-}) \, dM_t \\ & + (1 - H_t) G_t^{-1} (B_t \, dm_t - S_t \, d\mu_t) \\ & + (1 - H_t) G_t^{-2} (S_t \, d\langle\mu\rangle_t - B_t \, d\langle\mu, m\rangle_t), \end{aligned}$$

*where the $\mathbb{F}$-martingale $m$ is given by (15). Equivalently,*

(19)     $dS_t^c = r_t S_t^c \, dt + (Z_t - S_{t-}) \, dM_t + G_t^{-1} (B_t \, d\widehat{m}_t - S_t \, d\widehat{\mu}_t),$

*where the $\mathbb{G}$-martingales $\widehat{m}$ and $\widehat{\mu}$ are given by (12) and (13), respectively. The pre-default cumulative price $\widetilde{S}^c$ satisfies, for $t \in [0, T]$,*

(20)
$$\begin{aligned} d\widetilde{S}_t^c = {}& r_t \widetilde{S}_t^c \, dt + \lambda_t (\widetilde{S}_t - Z_t) \, dt \\ & + G_t^{-1} (B_t \, dm_t - \widetilde{S}_t \, d\mu_t) + G_t^{-2} (\widetilde{S}_t \, d\langle\mu\rangle_t - B_t \, d\langle\mu, m\rangle_t). \end{aligned}$$

PROOF.   Formula (10) yields

(21)
$$\begin{aligned} dS_t^c = {}& dS_t + d\Big(B_t \int_{]0,t]} B_u^{-1} \, dD_u\Big) = dS_t + r_t(S_t^c - S_t) \, dt + dD_t \\ = {}& dS_t + r_t(S_t^c - S_t) \, dt + (1 - H_t) \, dA_t + Z_t \, dH_t. \end{aligned}$$

By combining (22) with (17), we obtain (18). Formulae (19) and (20) are immediate consequences of (12), (13) and (18).  □

2.1.2.1. *Dynamics under Hypothesis (H).*   Let us now consider the special case where the so-called *Hypothesis*[1] $(H)$ is satisfied under $\mathbb{Q}^*$ between the filtrations $\mathbb{F}$ and $\mathbb{G} = \mathbb{H} \vee \mathbb{F}$. This means that the *immersion property* holds for the filtrations $\mathbb{F}$ and $\mathbb{G}$, in the sense that any $\mathbb{F}$-martingale under $\mathbb{Q}^*$ is also a $\mathbb{G}$-martingale under $\mathbb{Q}^*$. In that case, the survival process $G$ of $\tau$ with respect to $\mathbb{F}$ is known to be nonincreasing (see, e.g., Chapter 6 in [4] or [19]), so that $G = -\nu$. In other words, the continuous martingale $\mu$ in the Doob–Meyer decomposition of $G$ vanishes. Consequently, formula (14) becomes

(22)
$$\begin{aligned} dS_t = {}& -S_{t-} \, dM_t + (1 - H_t)((r_t S_t - \lambda_t Z_t) \, dt - dA_t) \\ & + (1 - H_t) B_t G_t^{-1} \, dm_t. \end{aligned}$$

---

[1]This property is referred to as the *martingale invariance property* of $\mathbb{F}$ and $\mathbb{G}$ in [4].

Similarly, (18) reduces to

$$(23) \qquad dS_t^c = r_t S_t^c \, dt + (Z_t - \widetilde{S}_t) \, dM_t + (1 - H_t) G_t^{-1} B_t \, dm_t$$

and (20) becomes

$$(24) \qquad d\widetilde{S}_t^c = r_t \widetilde{S}_t^c \, dt + \lambda_t (\widetilde{S}_t - Z_t) \, dt + G_t^{-1} B_t \, dm_t.$$

REMARK. Hypothesis (H) is a rather natural assumption in the present context. Indeed, it can be shown that it is necessarily satisfied under the postulate that the underlying $\mathbb{F}$-market model is complete and arbitrage-free, and the extended $\mathbb{G}$-market model is arbitrage-free (for details, see Blanchet–Scalliet and Jeanblanc [10]).

2.1.3. *Price dynamics of a CDS.* In Definition 2.5 of a stylized $T$-maturity credit default swap, we follow the convention adopted in [8]. Unlike in [8], the default protection stream is now represented by an $\mathbb{F}$-predictable process $\delta$. We assume that the default protection payment is received at the time of default and it equals $\delta_t$ if default occurs at time $t$ prior to or at maturity date $T$. Note that $\delta_t$ represents the protection payment, so that according to our notational convention, the recovery rate equals $1 - \delta_t$ rather than $\delta_t$. The notional amount of the CDS is equal to one monetary unit.

DEFINITION 2.5. The stylized $T$-maturity *credit default swap* (CDS) with a constant rate $\kappa$ and *recovery at default* is a defaultable claim $(0, A, Z, \tau)$ in which we set $Z_t = \delta_t$ and $A_t = -\kappa t$ for every $t \in [0, T]$. An $\mathbb{F}$-predictable process $\delta : [0, T] \to \mathbb{R}$ represents the *default protection* and a constant $\kappa$ is the fixed *CDS rate* (also termed the *spread* or *premium* of the CDS).

A credit default swap is thus a particular defaultable claim in which the promised payoff $X$ is null and the recovery process $Z$ is determined in reference to the estimated recovery rate of the reference credit name. We shall use the notation $D(\kappa, \delta, T, \tau)$ to denote the dividend process of a CDS. It follows immediately from Definition 2.2 that the dividend process $D(\kappa, \delta, T, \tau)$ of a stylized CDS equals, for every $t \in \mathbb{R}_+$,

$$
\begin{aligned}
(25) \qquad D_t(\kappa, \delta, T, \tau) &= \int_{]0, t \wedge T]} \delta_u \, dH_u - \kappa \int_{]0, t \wedge T]} (1 - H_u) \, du \\
&= \delta_\tau \mathbb{1}_{\{\tau \le t \wedge T\}} - \kappa (t \wedge T \wedge \tau).
\end{aligned}
$$

In a more realistic approach, the process $A$ is discontinuous, with jumps occurring at the premium payment dates. In this work we shall only deal with a stylized CDS with a continuously paid premium; for a more practical approach, we refer to Brigo [11] and Brigo and Morini [12].

Let us first examine the valuation formula for a stylized $T$-maturity CDS. Since we now have $X = 0, Z = \delta$ and $A_t = -\kappa t$, we deduce easily from (5) that the ex-dividend price (or *mark-to-market*) of such CDS contract equals, for every $t \in [0, T]$,

$$(26) \qquad S_t(\kappa, \delta, T, \tau) = \mathbb{1}_{\{t < \tau\}} (\widetilde{\delta}(t, T) - \kappa \widetilde{A}(t, T)),$$

where we denote, for any $t \in [0, T]$,

$$\widetilde{\delta}(t, T) = \frac{B_t}{G_t} \mathbb{E}_{\mathbb{Q}^*} (\mathbb{1}_{\{t < \tau \leq T\}} B_\tau^{-1} \delta_\tau \mid \mathcal{F}_t)$$

and

$$\widetilde{A}(t, T) = \frac{B_t}{G_t} \mathbb{E}_{\mathbb{Q}^*} \left( \int_t^{T \wedge \tau} B_u^{-1} du \; \Big| \; \mathcal{F}_t \right).$$

The quantity $\widetilde{\delta}(t, T)$ is the pre-default value at time $t$ of the protection leg, whereas $\widetilde{A}(t, T)$ represents the pre-default present value at time $t$ of one risky basis point paid up to the maturity $T$ or the default time $\tau$, whichever comes first. For ease of notation, we shall write $S_t(\kappa)$ in place of $S_t(\kappa, \delta, T, \tau)$ in what follows. Note that the quantities $\widetilde{\delta}(t, T)$ and $\widetilde{A}(t, T)$ are well defined at any date $t \in [0, T]$, and not only prior to default as the terminology "predefault values" might suggest.

We are in the position to state the following immediate corollary to Proposition 2.1.

COROLLARY 2.3.  *The ex-dividend price of a CDS equals, for any $t \in [0, T]$,*

$$(27) \qquad S_t(\kappa) = \mathbb{1}_{\{t < \tau\}} \frac{B_t}{G_t} \mathbb{E}_{\mathbb{Q}^*} \left( \int_t^T B_u^{-1} G_u (\delta_u \lambda_u - \kappa) \, du \; \Big| \; \mathcal{F}_t \right)$$

*and, thus, the cumulative price of a CDS equals, for any $t \in [0, T]$,*

$$(28) \qquad \begin{aligned} S_t^c(\kappa) &= \mathbb{1}_{\{t < \tau\}} \frac{B_t}{G_t} \mathbb{E}_{\mathbb{Q}^*} \left( \int_t^T B_u^{-1} G_u (\delta_u \lambda_u - \kappa) \, du \; \Big| \; \mathcal{F}_t \right) \\ &\quad + B_t \int_{]0, t]} B_u^{-1} \, dD_u. \end{aligned}$$

The next result is a direct consequence of Proposition 2.2 and Corollary 2.2.

COROLLARY 2.4.  *The dynamics of the ex-dividend price $S(\kappa)$ on $[0, T]$ are*

$$(29) \qquad \begin{aligned} dS_t(\kappa) = {}&- S_{t-}(\kappa) \, dM_t + (1 - H_t)(r_t S_t(\kappa) + \kappa - \lambda_t \delta_t) \, dt \\ &+ (1 - H_t) G_t^{-1} (B_t \, dn_t - S_t(\kappa) \, d\mu_t) \\ &+ (1 - H_t) G_t^{-2} (S_t(\kappa) \, d\langle \mu \rangle_t - B_t \, d\langle \mu, n \rangle_t), \end{aligned}$$

*where the $\mathbb{F}$-martingale $n$ is given by the formula*

$$(30) \qquad n_t = \mathbb{E}_{\mathbb{Q}^*}\left(\int_0^T B_u^{-1} G_u(\delta_u \lambda_u - \kappa)\, du \,\Big|\, \mathcal{F}_t\right).$$

*The cumulative price $S^c(\kappa)$ satisfies, for every $t \in [0, T]$,*

$$\begin{aligned} dS_t^c(\kappa) = {} & r_t S_t^c(\kappa)\, dt + (\delta_t - S_{t-}(\kappa))\, dM_t \\ & + (1 - H_t) G_t^{-1}(B_t\, dn_t - S_t(\kappa)\, d\mu_t) \\ & + (1 - H_t) G_t^{-2}(S_t(\kappa)\, d\langle\mu\rangle_t - B_t\, d\langle\mu, n\rangle_t), \end{aligned}$$

*or equivalently,*

$$(31) \quad dS_t^c(\kappa) = r_t S_t^c(\kappa)\, dt + (\delta_t - S_{t-}(\kappa))\, dM_t + G_t^{-1}(B_t\, d\widehat{n}_t - S_t(\kappa)\, d\widehat{\mu}_t),$$

*where the $\mathbb{G}$-martingales $\widehat{n}$ and $\widehat{\mu}$ are given by (12) and (13) respectively.*

**Dynamics under Hypothesis (H).** If the immersion property of $\mathbb{F}$ and $\mathbb{G}$ holds, the martingale $\mu$ is null and, thus, (29) reduces to

$$(32) \qquad \begin{aligned} dS_t(\kappa) = {} & -\widetilde{S}_t(\kappa)\, dM_t + (1 - H_t)(r_t S_t(\kappa) + \kappa - \lambda_t \delta_t)\, dt \\ & + (1 - H_t) B_t G_t^{-1}\, dn_t, \end{aligned}$$

since the process $\widetilde{S}_t(\kappa), t \in [0, T]$, is continuous and satisfies [cf. (16)]

$$(33) \qquad d\widetilde{S}_t(\kappa) = ((\lambda_t + r_t)\widetilde{S}_t(\kappa) + \kappa - \lambda_t \delta_t)\, dt + B_t G_t^{-1}\, dn_t.$$

Let us note that the quantity $\kappa - \lambda_t \delta_t$ can be informally interpreted as the *pre-default dividend rate* of a CDS.

Similarly, we obtain from (31)

$$(34) \qquad dS_t^c(\kappa) = r_t S_t^c(\kappa)\, dt + (\delta_t - \widetilde{S}_t(\kappa))\, dM_t + (1 - H_t) B_t G_t^{-1}\, dn_t$$

and

$$d\widetilde{S}_t^c(\kappa) = r_t \widetilde{S}_t^c(\kappa)\, dt + \lambda_t(\widetilde{S}_t(\kappa) - \delta_t)\, dt + B_t G_t^{-1}\, dn_t.$$

2.1.4. *Dynamics of the market CDS spread.* Let us now introduce the notion of the market CDS spread. It reflects the real-world feature that for any date $s$ the CDS issued at this time has the fixed spread chosen in such a way that the CDS is worthless at its inception. Note that the recovery process $\delta = (\delta_t)_{t \in [0,T]}$ is fixed throughout. We fix the maturity date $T$ and we assume that credit default swaps with different inception dates have a common recovery function $\delta$.

DEFINITION 2.6. The *$T$-maturity market CDS spread* $\kappa(s, T)$ at time $s \in [0, T]$ is the level of the CDS rate that makes the values of the two legs of a CDS equal to each other at time $s$.

It should be noted that CDSs are quoted in terms of spreads. At any date $t$, one can take at no cost a long or short position in the CDS issued at this date with the fixed rate equal to the actual value of the market CDS spread for a given maturity and a given reference credit name.

Let us stress that the market CDS spread $\kappa(s, T)$ is not defined neither at the moment of default nor after this date, so that we shall deal in fact with the pre-default value of the market CDS spread. Observe that $\kappa(s, T)$ is represented by an $\mathcal{F}_s$-measurable random variable. In fact, it follows immediately from (27) that $\kappa(s, T)$ admits the following representation, for any $s \in [0, T]$:

$$\kappa(s, T) = \frac{\widetilde{\delta}(s, T)}{\widetilde{A}(s, T)} = \frac{\mathbb{E}_{\mathbb{Q}^*}(\int_s^T B_u^{-1} G_u \delta_u \lambda_u \, du \mid \mathcal{F}_s)}{\mathbb{E}_{\mathbb{Q}^*}(\int_s^T B_u^{-1} G_u \, du \mid \mathcal{F}_s)} = \frac{K_s^1}{K_s^2},$$

where we denote

$$K_s^1 = \mathbb{E}_{\mathbb{Q}^*}\left(\int_s^T B_u^{-1} G_u \delta_u \lambda_u \, du \,\Big|\, \mathcal{F}_s\right)$$

and

$$K_s^2 = \mathbb{E}_{\mathbb{Q}^*}\left(\int_s^T B_u^{-1} G_u \, du \,\Big|\, \mathcal{F}_s\right).$$

In what follows, we shall write briefly $\kappa_s$ instead of $\kappa(s, T)$. The next result furnishes a convenient representation for the price at time $t$ of a CDS issued at some date $s \le t$, that is, the marked-to-market value of a CDS that exists already for some time (recall that the market value of the just issued CDS is null).

PROPOSITION 2.3. *The ex-dividend price $S(\kappa_s)$ of a $T$-maturity market CDS initiated at time $s$ equals, for every $t \in [s, T]$,*

$$(35) \qquad S_t(\kappa_s) = \mathbb{1}_{\{t < \tau\}}(\kappa_t - \kappa_s)\widetilde{A}(t, T) = \mathbb{1}_{\{t < \tau\}}\widetilde{S}_t(\kappa_s),$$

*where $\widetilde{S}_t(\kappa_s)$ is the pre-default ex-dividend price at time $t$.*

PROOF. To establish (35), it suffices to observe that $S_t(\kappa_s) = S_t(\kappa_s) - S_t(\kappa_t)$ since $S_t(\kappa_t) = 0$. Therefore, in order to conclude, it suffices to use (26) with $\kappa = \kappa_t$ and $\kappa = \kappa_s$. $\quad\square$

Let us now derive the dynamics of the market CDS spread. Let us define two $\mathbb{F}$-martingales

$$(36) \quad m_s^1 = \mathbb{E}_{\mathbb{Q}^*}\left(\int_0^T B_u^{-1} G_u \delta_u \lambda_u \, du \,\Big|\, \mathcal{F}_s\right) = K_s^1 + \int_0^s B_u^{-1} G_u \delta_u \lambda_u \, du$$

and

$$(37) \qquad m_s^2 = \mathbb{E}_{\mathbb{Q}^*} \left( \int_0^T B_u^{-1} G_u \, du \, \Big| \, \mathcal{F}_s \right) = K_s^2 + \int_0^s B_u^{-1} G_u \, du.$$

Under Assumption 2.2, the $\mathbb{F}$-martingales $m^1$ and $m^2$ are continuous. Therefore, using the Itô formula, we find easily that the semimartingale decomposition of the market spread process reads

$$d\kappa_s = \frac{1}{K_s^2} \Big( B_s^{-1} G_s (\kappa_s - \delta_s \lambda_s) \, ds + \frac{\kappa_s}{K_s^2} \, d\langle m^2 \rangle_s - \frac{1}{K_s^2} \, d\langle m^1, m^2 \rangle_s \Big)$$

$$+ \frac{1}{K_s^2} \big( dm_s^1 - \kappa_s \, dm_s^2 \big).$$

2.2. *Replication of a defaultable claim.* We now assume that $k$ credit default swaps with certain maturities $T_i \geq T$ spreads $\kappa_i$ and protection payments $\delta^i$ for $i = 1, \ldots, k$ are traded over the time interval $[0, T]$. All these contracts are supposed to refer to the same underlying credit name and, thus, they have a common default time $\tau$. Formally, this family of CDSs is represented by the associated dividend processes $D^i = D(\kappa_i, \delta^i, T_i, \tau)$ given by formula (25). For brevity, the corresponding ex-dividend price will be denoted as $S^i(\kappa_i)$ rather than $S(\kappa_i, \delta^i, T_i, \tau)$. Similarly, $S^{c,i}(\kappa_i)$ will stand for the cumulative price process of the $i$th traded CDS. The 0th traded asset is the savings account $B$.

2.2.1. *Self-financing trading strategies in the CDS market.* Our goal is to examine hedging strategies for a defaultable claim $(X, A, Z, \tau)$. As expected, we will trade in $k$ credit default swaps and the savings account. To this end, we will consider trading strategies $\varphi = (\varphi^0, \ldots, \varphi^k)$, where $\varphi^0$ is a $\mathbb{G}$-adapted process and the processes $\varphi^1, \ldots, \varphi^k$ are $\mathbb{G}$-predictable.

In the present set-up we consider trading strategies that are self-financing in the standard sense, as recalled in the following definition.

DEFINITION 2.7. The *wealth process* $V(\varphi)$ of a strategy $\varphi = (\varphi^0, \ldots, \varphi^k)$ in the savings account $B$ and ex-dividend CDS prices $S^i(\kappa_i), i = 1, \ldots, k$ equals, for any $t \in [0, T]$,

$$(38) \qquad V_t(\varphi) = \varphi_t^0 B_t + \sum_{i=1}^k \varphi_t^i S_t^i(\kappa_i).$$

A strategy $\varphi$ is said to be *self-financing* if $V_t(\varphi) = V_0(\varphi) + G_t(\varphi)$ for every $t \in [0, T]$, where the gains process $G(\varphi)$ is defined as follows:

$$(39) \qquad G_t(\varphi) = \int_{]0,t]} \varphi_u^0 \, dB_u + \sum_{i=1}^k \int_{]0,t]} \varphi_u^i \, d(S_u^i(\kappa_i) + D_u^i),$$

where $D^i = D(\kappa_i, \delta^i, T_i, \tau)$ is the dividend process of the $i$th CDS [see formula (25)].

The following lemma is fairly general; in particular, it is independent of the choice of the underlying model. Indeed, in the proof of this result we only use the obvious relationships $dB_t = r_t B_t\, dt$ and the relationship [cf. (10)]

$$(40) \qquad S_t^{c,i}(\kappa_i) = S_t^i(\kappa_i) + B_t \int_{]0,t]} B_u^{-1}\, dD_u^i.$$

Let $V^*(\varphi) = B^{-1} V(\varphi)$ be the discounted wealth process and let $S^{c,i,*}(\kappa_i) = B^{-1} S^{c,i}(\kappa_i)$ be the discounted cumulative price.

LEMMA 2.2. *Let* $\varphi = (\varphi^0, \ldots, \varphi^k)$ *be a self-financing trading strategy in the savings account* $B$ *and ex-dividend prices* $S^i(\kappa_i), i = 1, \ldots, k$. *Then the discounted wealth process* $V^* = B^{-1} V(\varphi)$ *satisfies, for* $t \in [0, T]$,

$$(41) \qquad dV_t^*(\varphi) = \sum_{i=1}^k \varphi_t^i\, dS_t^{c,i,*}(\kappa_i).$$

PROOF.   We have

$$dV_t^*(\varphi) = B_t^{-1}\, dV_t(\varphi) - r_t B_t^{-1} V_t(\varphi)\, dt = B_t^{-1}(dV_t(\varphi) - r_t V_t(\varphi)\, dt)$$

$$= B_t^{-1} \left[ \varphi_t^0 r_t B_t\, dt + \sum_{i=1}^k \varphi_t^i (dS_t^i(\kappa_i) + dD_t^i) - r_t V_t(\varphi)\, dt \right]$$

$$= B_t^{-1} \left[ \left( V_t(\varphi) - \sum_{i=1}^k \varphi_t^i S_t^i(\kappa_i) \right) r_t\, dt \right.$$

$$\qquad\qquad \left. + \sum_{i=1}^k \varphi_t^i (dS_t^i(\kappa_i) + dD_t^i) - r_t V_t(\varphi)\, dt \right]$$

$$= B_t^{-1} \sum_{i=1}^k \varphi_t^i (dS_t^i(\kappa_i) - r_t S_t^i(\kappa_i)\, dt + dD_t^i)$$

$$= \sum_{i=1}^k \varphi_t^i (d(B_t^{-1} S_t^i(\kappa_i)) + B_t^{-1}\, dD_t^i).$$

By comparing the last formula with (40), we conclude that (41) holds.   □

2.2.2. *Replication with ex-dividend prices of CDSs.*   Recall that the cumulative price of a defaultable claim $(X, A, Z, \tau)$ is denoted as $S^c$. We adopt the following natural definition of replication of a defaultable claim. Note that the set of traded assets is not explicitly specified in this definition. Hence, it can be used for any choice of primary traded assets.

DEFINITION 2.8. We say that a self-financing strategy $\varphi = (\varphi^0, \ldots, \varphi^k)$ *replicates* a defaultable claim $(X, A, Z, \tau)$ if its wealth process $V(\varphi)$ satisfies $V_t(\varphi) = S_t^c$ for every $t \in [0, T]$. In particular, the equality $V_{t \wedge \tau}(\varphi) = S_{t \wedge \tau}^c$ holds for every $t \in [0, T]$.

In the remaining part of this section we assume that Hypothesis (H) holds. Hence, the hazard process $\Gamma$ is increasing and, thus, by Assumption 2.1, we have that, for any $t \in [0, T]$,

$$\Gamma_t = \Lambda_t = \int_0^t \lambda_u \, du.$$

The discounted cumulative price $S^{c,i}(\kappa_i)$ of the $i$th CDS is governed by [cf. (34)]

$$(42) \qquad dS_t^{c,i,*}(\kappa_i) = B_t^{-1}(\delta_t^i - \widetilde{S}_t^i(\kappa_i)) \, dM_t + (1 - H_t)G_t^{-1} \, dn_t^i,$$

where [cf. (30)]

$$(43) \qquad n_t^i = \mathbb{E}_{\mathbb{Q}^*}\left( \int_0^{T_i} B_u^{-1} G_u(\delta_u^i \lambda_u - \kappa_i) \, du \,\Big|\, \mathcal{F}_t \right).$$

The next lemma yields the dynamics of the wealth process $V(\varphi)$ for a self-financing strategy $\varphi$.

LEMMA 2.3. *For any self-financing trading strategy $\varphi$, the discounted wealth $V^*(\varphi) = B^{-1} V(\varphi)$ satisfies, for any $t \in [0, T]$,*

$$(44) \qquad dV_t^*(\varphi) = \sum_{i=1}^k \varphi_t^i (B_t^{-1}(\delta_t^i - \widetilde{S}_t^i(\kappa_i)) \, dM_t + (1 - H_t)G_t^{-1} \, dn_t^i).$$

PROOF. It suffices to combine (41) with (42). □

It is clear from the lemma that it is enough to search for the components $\varphi^1, \ldots, \varphi^k$ of a strategy $\varphi$. The same remark applies to self-financing strategies introduced in Definitions 2.7 and 2.10 below.

It is worth stressing that in what follows we shall only consider *admissible* trading strategies, that is, strategies for which the discounted wealth process $V^*(\varphi) = B^{-1} V(\varphi)$ is a $\mathbb{G}$-martingale under $\mathbb{Q}^*$. The market model in which only admissible trading strategies are allowed is *arbitrage-free*, that is, arbitrage opportunities are ruled out. Admissibility of a replicating strategy will be ensured by the equality $V(\varphi) = S^c$ and the fact that the discounted cumulative price $B^{-1} S^c$ of a defaultable claim is a $\mathbb{G}$-martingale under $\mathbb{Q}^*$.

We work throughout under the standing Assumptions 2.1 and 2.2 and the following postulate.

ASSUMPTION 2.3. *The filtration $\mathbb{F}$ is generated by a $d$-dimensional Brownian motion $W$ under $\mathbb{Q}^*$.*

Since Hypothesis (H) is assumed to hold, the process $W$ is also a Brownian motion with respect to the enlarged filtration $\mathbb{G} = \mathbb{H} \vee \mathbb{F}$. Recall that all (local) martingales with respect to a Brownian filtration are necessarily continuous. Hence, Assumption 2.2 is obviously satisfied.

The crucial observation is that, by the predictable representation property of a Brownian motion, there exist $\mathbb{F}$-predictable, $\mathbb{R}^d$-valued processes $\xi$ and $\zeta^i, i = 1, \ldots, k$, such that $dm_t = \xi_t \, dW_t$ and $dn_t^i = \zeta_t^i \, dW_t$, where $m$ and $n^i$ are given by (15) and (43) respectively.

We are now in the position to state the hedging result for a defaultable claim in the single-name set-up. We consider a defaultable claim $(X, A, Z, \tau)$ satisfying the natural integrability conditions under $\mathbb{Q}^*$, such that the cumulative price process $S^c$ for this claim is well defined.

THEOREM 2.1. *Assume that there exist $\mathbb{F}$-predictable processes $\varphi^1, \ldots, \varphi^k$ satisfying the following conditions, for any $t \in [0, T]$:*

$$\text{(45)} \qquad \sum_{i=1}^{k} \varphi_t^i (\delta_t^i - \widetilde{S}_t^i(\kappa_i)) = Z_t - \widetilde{S}_t, \qquad \sum_{i=1}^{k} \varphi_t^i \zeta_t^i = \xi_t.$$

*Let the process $V(\varphi)$ be given by (44) with the initial condition $V_0(\varphi) = S_0^c$ and let $\varphi^0$ be given by, for $t \in [0, T]$,*

$$\text{(46)} \qquad \varphi_t^0 = B_t^{-1} \left( V_t(\varphi) - \sum_{i=1}^{k} \varphi_t^i S_t^i(\kappa_i) \right).$$

*Then the self-financing trading strategy $\varphi = (\varphi^0, \ldots, \varphi^k)$ in the savings account $B$ and assets $S^i(\kappa_i), i = 1, \ldots, k$, replicates the defaultable claim $(X, A, Z, \tau)$.*

PROOF. From Lemma 2.3, we know that the discounted wealth process satisfies

$$\text{(47)} \qquad dV_t^*(\varphi) = \sum_{i=1}^{k} \varphi_t^i (B_t^{-1}(\delta_t^i - \widetilde{S}_t^i(\kappa_i)) \, dM_t + (1 - H_t) G_t^{-1} \, dn_t^i).$$

Recall also that the discounted cumulative price $S^{c,*}$ of a defaultable claim is governed by [cf. (23)]

$$\text{(48)} \qquad dS_t^{c,*} = B_t^{-1}(Z_t - \widetilde{S}_t) \, dM_t + (1 - H_t) G_t^{-1} \, dm_t.$$

We will show that if the two conditions in (45) are satisfied for any $t \in [0, T]$, then the equality $V_t(\varphi) = S_t^c$ holds for any $t \in [0, T]$.

Let $\widetilde{V}^*(\varphi) = B^{-1}\widetilde{V}(\varphi)$ stand for the *discounted pre-default wealth*, where $\widetilde{V}(\varphi)$ is the unique $\mathbb{F}$-adapted process such that $\mathbb{1}_{\{t<\tau\}}V_t(\varphi) = \mathbb{1}_{\{t<\tau\}}\widetilde{V}_t(\varphi)$ for every $t \in [0,T]$. On the one hand, using (45), we obtain

$$d\widetilde{V}_t^*(\varphi) = \sum_{i=1}^{k} \varphi_t^i (\lambda_t B_t^{-1}(\widetilde{S}_t^i(\kappa_i) - \delta_t^i)\,dt + G_t^{-1}\zeta_t^i\,dW_t)$$

$$= \lambda_t B_t^{-1}(\widetilde{S}_t - Z_t)\,dt + G_t^{-1}\xi_t\,dW_t.$$

On the other hand, the discounted pre-default cumulative price $\widetilde{S}^{c,*} = B^{-1}\widetilde{S}^c$ satisfies [cf. (24)]

$$d\widetilde{S}_t^{c,*} = \lambda_t B_t^{-1}(\widetilde{S}_t - Z_t)\,dt + G_t^{-1}\xi_t\,dW_t.$$

Since by assumption $\widetilde{V}_0^*(\varphi) = V_0(\varphi) = S_0^c = \widetilde{S}_0^{c,*}$, it is clear that $\widetilde{V}_t^*(\varphi) = \widetilde{S}_t^{c,*}$ for every $t \in [0,T]$. We thus conclude that the pre-default wealth $\widetilde{V}(\varphi)$ of $\varphi$ and the pre-default cumulative price $\widetilde{S}^c$ of the claim coincide. Note that the first equality in (45) is in fact only essential for those values of $t \in [0,T]$ for which $\lambda_t \neq 0$.

To complete the proof, we need to check what happens when default occurs prior to or at maturity $T$. To this end, it suffices to compare the jumps of $S^c$ and $V(\varphi)$ at time $\tau$. In view of (47), (48) and (45), we obtain

$$\Delta V_\tau(\varphi) = Z_\tau - \widetilde{S}_\tau = \Delta S_\tau^c$$

and, thus, $V_{t\wedge\tau}(\varphi) = S_{t\wedge\tau}^c$ for any $t \in [0,T]$. After default, we have $dV_t(\varphi) = r_t V_t(\varphi)\,dt$ and $dS_t^c = r_t S_t^c\,dt$, so that we conclude that the desired equality $V_t(\varphi) = S_t^c$ holds for any $t \in [0,T]$.  $\square$

As pointed out by the anonymous referee, a BSDE approach can be useful for the determination of the quantities $\zeta^i$ and $\xi$. The coefficient of the jumping martingale in the BSDE is determined explicitly by the size of the jump, as noticed in the general representation theorem presented in Section 6.3. of Bielecki and Rutkowski [4].

2.2.3. *Replication with rolling CDSs.* When considering trading strategies involving CDSs that were issued in the past, one encounters a practical difficulty regarding their liquidity. For this reason, we shall now analyze trading strategies based on rolling CDS contracts. Toward this end, we will define a contract—that we call a *rolling credit default swap*—which at any time $t$ has similar features as the $T$-maturity CDS issued at this date $t$, in particular, its ex-dividend price is equal to zero. Intuitively, one can think of the rolling CDS of a constant maturity $T$ as a stream of CDSs of constant maturities equal to $T$ that are continuously entered into and immediately unwound. Thus, a rolling CDS contract is equivalent to a self-financing trading strategy that at any given time $t$ enters into a CDS contract of maturity $T$ and then unwinds the contract at time $t + dt$.

REMARK. We maintain the assumptions of Section 2.2.2. Also, we use here a simplifying assumption that the protection payment process $\delta$ is generic, that is, it is the same for all CDS contracts referencing the same default $\tau$ and with the same maturity $T$. Otherwise, for every fixed maturity date $T$, we would need to consider the whole class of protection payment processes, indexed by the initiation date.

Denote by $S_u(t, \kappa(t, T))$ the time—$u$ ex-dividend price of a CDS that was initiated at time $t \le u$ at the contracted spread $\kappa(t, T)$, where $T$ is the maturity. From (32), we immediately obtain that, for any fixed $t \in [0, T]$ and every $u \in [t, T]$,

$$
\begin{aligned}
dS_u(t, \kappa(t, T)) = &-\widetilde{S}_u(t, \kappa(t, T))\, dM_u \\
&+ (1 - H_u)(r_u \widetilde{S}_u(t, \kappa(t, T)) + \kappa(t, T) - \lambda_u \delta_u)\, du \\
&+ (1 - H_u) B_u G_u^{-1}(dn_u^1(t) - \kappa(t, T)\, dn_u^2(t)),
\end{aligned}
$$

where we denote

$$
n_u^1(t) = \mathbb{E}_{\mathbb{Q}^*}\left( \int_t^T B_s^{-1} G_s \delta_s \lambda_s \, ds \; \Big| \; \mathcal{F}_u \right)
$$

and

$$
n_u^2(t) = \mathbb{E}_{\mathbb{Q}^*}\left( \int_t^T B_s^{-1} G_s \, ds \; \Big| \; \mathcal{F}_u \right).
$$

Note also that for fixed $t \in [0, T]$ and every $u > t$ [see (36) and (37)],

$$
\begin{aligned}
n_u^1(t) &= \mathbb{E}_{\mathbb{Q}^*}\left( \int_0^T B_s^{-1} G_s \delta_s \lambda_s \, ds \; \Big| \; \mathcal{F}_u \right) - \int_0^t B_s^{-1} G_s \delta_s \lambda_s \, ds \\
&= m_u^1 - \int_0^t B_s^{-1} G_s \delta_s \lambda_s \, ds,
\end{aligned}
$$

so that $dn_u^1(t) = dm_u^1$. A similar argument shows that $dn_u^2(t) = dm_u^2$. Consequently, we also have that

$$
\begin{aligned}
dS_u(t, \kappa(t, T)) = &-\widetilde{S}_u(t, \kappa(t, T))\, dM_u \\
&+ (1 - H_u)(r_u S_u(t, \kappa(t, T)) + \kappa(t, T) - \lambda_u \delta_u)\, du \\
&+ (1 - H_u) B_u G_u^{-1}(dm_u^1 - \kappa(t, T)\, dm_u^2).
\end{aligned}
\tag{49}
$$

We are in the position to prove the following result, in which $R$ stands for the wealth process of a self-financing trading strategy representing the rolling CDS with maturity date $T$. Note, in particular, that the process $B_t^{-1} R_t$ is a local martingale, in general, and a proper martingale under mild assumptions, which we take for granted.

LEMMA 2.4.   *The dynamics of the wealth process $R$ are*

$$(50) \qquad dR_t = r_t R_t \, dt + \delta_t \, dM_t + (1 - H_t) B_t G_t^{-1} (dm_t^1 - \kappa(t, T) \, dm_t^2).$$

PROOF.   Let us sketch the proof of the lemma. Assume that we enter in a CDS at time $t$ and we close this position at time $t + h$, putting all the remaining money in the bank account and entering in a new CDS, issued at time $t + h$. Since the savings account clearly accrues at the short-term rate $r$, we may postulate, for simplicity of presentation, that $R_t = 0$. Also, it is enough to examine the dynamics of the wealth process $R$ up to the moment $\tau \wedge T$. For instance, the wealth of our portfolio at time $s = (t + h) \wedge \tau$ can be represented as follows:

$$R_s = S_s(t, \kappa(t, T)) - \kappa(t, T) B_s \int_t^s B_u^{-1} \, du + \delta_\tau \mathbb{1}_{\{t < \tau \le s\}}.$$

For $\Delta R = R_s - R_t = R_s$, we thus obtain [recall that $S_t(t, \kappa(t, T)) = 0$]

$$\begin{aligned}
\Delta R &= -\int_t^s \widetilde{S}_u(t, \kappa(t, T)) \, dM_u + \int_t^s (r_u \widetilde{S}_u(t, \kappa(t, T)) + \kappa(t, T) - \lambda_u \delta_u) \, du \\
&\quad + \int_t^s B_u G_u^{-1} (dm_u^1 - \kappa(t, T) \, dm_u^2) - \int_t^s B_s B_u^{-1} \kappa(t, T) \, du + \int_t^s \delta_u \, dH_u \\
&= \int_t^s r_u \widetilde{S}_u(t, \kappa(t, T)) \, du + \int_t^s (\delta_u - \widetilde{S}_u(t, \kappa(t, T))) \, dM_u \\
&\quad + \int_t^s (1 - B_s B_u^{-1}) \kappa(t, T) \, du + \int_t^s B_u G_u^{-1} (dm_u^1 - \kappa(t, T) \, dm_u^2).
\end{aligned}$$

Recall that under our assumptions the process $\kappa(t, T), t \in [0, T]$ is $\mathbb{F}$-adapted and continuous, and thus $\mathbb{F}$-predictable (see Section 2.1.4). Hence, by letting $h$ tend to zero, we see that the process $R$ stopped at $\tau$ satisfies (recall that we postulated that $R_t = 0$)

$$dR_t = \delta_t \, dM_t + B_t G_t^{-1} (dm_t^1 - \kappa(t, T) \, dm_t^2).$$

In general, we obtain the asserted formula (50).   $\square$

DEFINITION 2.9.   Let $0 < T_1 < T_2$. The *rolling CDS* with initial time $T_1$ and expiry date $T_2$ is a financial security initiated at time $T_1$, whose ex-dividend price is zero and whose cumulative dividend process, say, $R^c(T_1, T_2)$, satisfies $R_{T_1}^c(T_1, T_2) = 0$ and, for every $t \in [T_1, T_2]$,

$$\begin{aligned}
(51) \qquad dR_t^c(T_1, T_2) &= r_t R_t^c(T_1, T_2) \, dt + \delta_t \, dM_t \\
&\quad + (1 - H_t) B_t G_t^{-1} (dm_t^1 - \kappa(t, T_2) \, dm_t^2).
\end{aligned}$$

In view of Lemma 2.4, the process $R(T_1, T_2)$ represents the wealth of a strategy in which, at any time $t \in [T_1, T_2]$ prior to default time $\tau$, we enter into one long CDS contract initiated at time $t$ and maturing at time $T_2$ and then immediately (i.e., at time $t + dt$) we unwind this contract and we enter into one long CDS contract initiated at time $t + dt$ and maturing at time $T_2$. As we shall now argue, practical considerations suggest that the actual life-span of a rolling CDS used for hedging purposes should be less than its nominal duration $[T_1, T_2]$.

In market practice the tenor of maturities of CDS contracts consists of, typically, four dates per calendar year. This means that, for example, 3-year CDSs issued between October 1 and December 31, 2007, will all mature on December 31, 2010; likewise, 3-year CDS contracts issued between January 1 and March 31, 2008, will all mature on March 31, 2011. Consequently, we need to deal with a sequence of rolling CDSs, all with the life-span $\mathcal{T}$ of 3 months. For instance, a rolling CDS started on $T_1^k = $ October 1, 2007 and maturing on $T_2^k = $ December, 31, 2010 will need to be replaced on January 1, 2008, by a new rolling CDS starting on $T_1^{k+1} = $ January 1, 2008 and maturing on $T_2^{k+1} = $ March, 31, 2011, and so on. Analogous rules apply to 5-year, 7-year and 10-year credit default swaps.

We shall now introduce a portfolio of rolling CDSs, all with the same life-span $\mathcal{T}$ of 3 months. We shall denote the contracts in the portfolio by $R^{k,i}$ with spreads $\kappa^{k,i}$ and protection payments $\delta^{k,i}$, where $k = 0, \ldots, K$ and $i = 1, \ldots, I_k$. The superscript $k$ determines the initial date of the given rolling CDS, say, $T_1^{k,i}$; we assume that $T_1^{k,i} = k\mathcal{T}$ for all $i = 1, \ldots, I_k$. The superscript $i$ identifies other relevant features of the contract, such as the nominal duration of the contract, for example, three years, five years, etc. For example, we may have that $T_1^{k,1} = T_1^{k,2} = $ October 1, 2007 and $T_2^{k,1} = $ December 31, 2010 (a class of three year contracts) and $T_2^{k,2} = $ December 31, 2012 (a class of five year contracts). Note that the contract $R^{k,i}$ is only alive for $t \in [k\mathcal{T}, (k+1)\mathcal{T})$ and, thus, we formally set $R^{k,i}$ to be zero for $t \notin [k\mathcal{T}, (k+1)\mathcal{T})$.

To describe the self-financing trading strategies in the savings account $B$ and rolling CDS contracts, we may use the following version of Definition 2.7.

DEFINITION 2.10.   A strategy $\varphi = (\varphi^0, \{\varphi^{k,i}, k = 0, \ldots, K, i = 1, \ldots, I_k\})$ in the savings account $B$ and the rolling CDS contracts with the cumulative dividend processes $D^{k,i} = R^{k,i}, k = 0, \ldots, K, i = 1, \ldots, I_k$, is said to be *self-financing* if the wealth $V_t(\varphi) = \varphi_t^0 B_t$ satisfies $V_t(\varphi) = V_0(\varphi) + G_t(\varphi)$ for every $t \in [0, T]$, where the gains process $G(\varphi)$ is defined as follows:

$$G_t(\varphi) = \int_{]0,t]} \varphi_u^0 \, dB_u + \sum_{k,i} \int_{]0,t]} \varphi_u^{k,i} \, dR_u^{k,i}.$$

Similarly as in Lemma 2.2, we obtain the following condition satisfied by the discounted wealth process of any self-financing strategy $\varphi$ in the sense of Definition 2.10:

$$(52) \qquad d(B_t^{-1} V_t(\varphi)) = \sum_{k,i=1} \varphi_t^{k,i} \, d(B_t^{-1} R_t^{k,i}).$$

Let us denote by $m^{1,k,i}$ the process given by formula (36) with $\delta$ replaced by $\delta^{k,i}$. Using (51), we obtain, on $[0, \tau \wedge T]$,

$$(53) \qquad \begin{aligned} d(B_t^{-1} R_t^{k,i}) &= G_t^{-1}(dm_t^{1,k,i} - \kappa^{k,i}(t)\, dm_t^2) + B_t^{-1}\delta_t^{k,i}\, dM_t \\ &= G_t^{-1}(\zeta_t^{1,k,i} - \kappa^{k,i}(t)\zeta_t^2)\, dW_t + B_t^{-1}\delta_t^{k,i}\, dM_t, \end{aligned}$$

for some $\mathbb{F}$-predictable processes $\zeta^{1,k,i}$ and $\zeta^2$. The following results can be demonstrated in a way analogous to the proof of Theorem 2.1 and, thus, the details are left to the reader.

THEOREM 2.2. *Assume that there exist $\mathbb{F}$-predictable processes $\varphi^{k,i}$ satisfying the following conditions, for every $t \in [0,T]$:*

$$(54) \qquad \begin{aligned} \sum_{k,i} \mathbb{1}_{[k\mathcal{T},(k+1)\mathcal{T})}(t)\varphi_t^{k,i}\delta_t^{k,i} &= Z_t - \widetilde{S}_t, \\ \sum_{k,i} \varphi_t^{k,i}(\zeta_t^{1,k,i} - \kappa^{k,i}(t)\zeta_t^2) &= \xi_t. \end{aligned}$$

*Let the process $V(\varphi)$ be given by (52) with the initial condition $V_0(\varphi) = S_0^c$ and let $\varphi^0$ satisfy $\varphi_t^0 B_t = V_t(\varphi)$ for every $t \in [0,T]$. Then the self-financing trading strategy $\varphi$ in the savings account $B$ and rolling credit default swaps $R^{k,i}, k = 0, \ldots, K, i = 1, \ldots, I_k$, replicates the defaultable claim $(X, A, Z, \tau)$.*

REMARK. Laurent, Cousin and Fermanian [22] take a different approach to hedging with market CDSs. They postulate the existence of *instantaneous digital CDSs*. Formally, the instantaneous dynamics of such a contract is exactly the same as these of the martingale $M$. Note, however, that they deal with a model in which the reference filtration is trivial and this corresponds to the set-up considered in our previous paper [8].

2.2.4. *Sufficient conditions for hedgeability.* The first equality in (45) eliminates the *jump risk*, whereas the second one is used to eliminate the *spread risk*. In general, the existence of $\varphi^1, \ldots, \varphi^k$ satisfying (45) is not ensured and it is easy to give an example when a solution to (45) fails to exist. In Example 2.1 below, we deal with a (admittedly somewhat artificial) situation when the jump risk can be perfectly hedged, but the prices of traded

CDSs are deterministic prior to default, so that the spread risk of a default-able claim is nonhedgeable. In general, the solvability of (45) depends on several factors, such as: the number of traded assets, the dimension of the driving Brownian motion, the random character of default intensity $\gamma$ and recovery payoffs $\delta^i$ and the features of a defaultable claim that we wish to hedge.

EXAMPLE 2.1. Let $r = 0$ and $k = 2$. Assume that $\kappa_1 \neq \kappa_2$ are nonzero constants, $T_1 \neq T_2$, and let $\delta^1 = \delta^2 = Z = 0$. Assume also that the default intensity $\gamma(t) > 0$ is deterministic and the promised payoff $X$ is a nonconstant $\mathcal{F}_T$-measurable random variable. We thus have [cf. (15)]

$$m_t = \mathbb{E}_{\mathbb{Q}^*}(G_T X \mid \mathcal{F}_t) = G_T \mathbb{E}_{\mathbb{Q}^*}(X \mid \mathcal{F}_t) = G_T\left(\mathbb{E}_{\mathbb{Q}^*}(X) + \int_0^t \xi_u \, dW_u\right)$$

for some nonvanishing process $\xi$. However, since $\gamma$ is deterministic, it is also easy to deduce from (30) that $\zeta_t^i = 0, i = 1, 2$, for every $t \in [0, T]$. The first condition in (45) reads $\sum_{i=1}^2 \varphi_t^i \widetilde{S}_t^i(\kappa_i) = \widetilde{S}_t$, and since manifestly $\widetilde{S}_t^i(\kappa_i) = \kappa_i G_t \int_t^T G_u \, du \neq 0$ for every $t \in [0, T]$, no difficulty may arise here. However, the second equality $\sum_{i=1}^2 \varphi_t^i \zeta_t^i = \xi_t$ cannot be satisfied for every $t \in [0, T]$, since the left-hand side vanishes for every $t \in [0, T]$.

We shall now provide sufficient conditions for the existence and uniqueness of a replicating strategy for any defaultable claim in the practically appealing case of CDSs with constant protection payments. We first address this issue in the special case where $k = 2$ and the model is driven by a one-dimensional Brownian motion $W$. In addition, we assume that the two traded CDSs have the same maturity, $T_1 = T_2 = U$; this assumption is made here for simplicity of presentation only, and it will be relaxed in Proposition 2.5 below. Let us denote

$$\widetilde{P}_t = B_t G_t^{-1} \mathbb{E}_{\mathbb{Q}^*}\left(\int_t^U B_u^{-1} G_u \lambda_u \, du \, \Big| \, \mathcal{F}_t\right),$$

$$\widetilde{A}_t = B_t G_t^{-1} \mathbb{E}_{\mathbb{Q}^*}\left(\int_t^U B_u^{-1} G_u \, du \, \Big| \, \mathcal{F}_t\right)$$

so that $\widetilde{S}_t^i(\kappa_i) = \delta_i \widetilde{P}_t - \kappa_i \widetilde{A}_t$ for $i = 1, 2$. Similarly [cf. (43)], $n_t^i = \delta_i m_t^1 - \kappa_i m_t^2$, where we set

$$m_t^1 = \mathbb{E}_{\mathbb{Q}^*}\left(\int_0^U B_u^{-1} G_u \lambda_u \, du \, \Big| \, \mathcal{F}_t\right), \qquad m_t^2 = \mathbb{E}_{\mathbb{Q}^*}\left(\int_0^U B_u^{-1} G_u \, du \, \Big| \, \mathcal{F}_t\right).$$

By the predictable representation property of the Brownian motion, $dm_t^j = \psi_t^j \, dW_t$ for $j = 1, 2$ for some $\mathbb{F}$-predictable, real-valued processes $\psi^1$ and $\psi^2$.

PROPOSITION 2.4. *Assume that $T_1 = T_2 = U$ and the constant protection payments $\delta_1$ and $\delta_2$ are such that $\delta_1 \kappa_2 \neq \delta_2 \kappa_1$ and $(1 - \tilde{P}_t) \psi_t^2 \neq \tilde{A}_t \psi_t^1$ for almost every $t \in [0, T]$. Then for any defaultable claim $(X, A, Z, \tau)$ there exists a unique solution $(\varphi^1, \varphi^2)$ to (45).*

PROOF. In view of (43), we obtain $\zeta_t^i = \delta_i \psi_t^1 - \kappa_i \psi_t^2$. Hence, the matching conditions (45) become

$$(55) \qquad \sum_{i=1}^{2} \varphi_t^i (\delta_i \hat{P}_t + \kappa_i \tilde{A}_t) = Z_t - \tilde{S}_t, \qquad \sum_{i=1}^{2} \varphi_t^i (\delta_i \psi_t^1 - \kappa_i \psi_t^2) = \xi_t,$$

where, for brevity, we write $\hat{P}_t = 1 - \tilde{P}_t$. A unique solution to (55) exists provided that the random matrix

$$N_t = \begin{bmatrix} \delta_1 \hat{P}_t + \kappa_1 \tilde{A}_t & \delta_2 \hat{P}_t + \kappa_2 \tilde{A}_t \\ \delta_1 \psi_t^1 - \kappa_1 \psi_t^2 & \delta_2 \psi_t^1 - \kappa_2 \psi_t^2 \end{bmatrix}$$

is nonsingular for almost every $t \in [0, T]$, that is, whenever

$$\det N_t = (\delta_2 \kappa_1 - \delta_1 \kappa_2)(\hat{P}_t \psi_t^2 - \tilde{A}_t \psi_t^1) \neq 0$$

for almost every $t \in [0, T]$.  □

Equality $\delta_2 \kappa_1 - \delta_1 \kappa_2 = 0$ would practically mean that we deal with a single CDS rather than two distinct CDSs. Note that we have here two sources of uncertainty, the discontinuous martingale $M$ and the Brownian motion $W$; hence, it was natural to expect that the number of assets required to span the market equals 3.

If the model is driven by a $d$-dimensional Brownian motion, it is natural to expect that one will need at least $d + 2$ assets (the savings account and $d + 1$ distinct CDSs, say) to replicate any defaultable claim, that is, to ensure the model's completeness (for similar results in a Markovian set-up, see [23]). This question is examined in the next result, in which we denote

$$\widetilde{P}_t^i = B_t G_t^{-1} \mathbb{E}_{\mathbb{Q}^*} \left( \int_t^{T_i} B_u^{-1} G_u \lambda_u \, du \,\Big|\, \mathcal{F}_t \right),$$

$$\widetilde{A}_t^i = B_t G_t^{-1} \mathbb{E}_{\mathbb{Q}^*} \left( \int_t^{T_i} B_u^{-1} G_u \, du \,\Big|\, \mathcal{F}_t \right),$$

so that $\widetilde{S}_t^i(\kappa_i) = \delta_i \widetilde{P}_t^i - \kappa_i \widetilde{A}_t^i$ for $i = 1, \ldots, k$. Similarly [cf. (43)], $n_t^i = \delta_i m_t^{1i} - \kappa_i m_t^{2i}$, where we set

$$m_t^{1i} = \mathbb{E}_{\mathbb{Q}^*} \left( \int_0^{T_i} B_u^{-1} G_u \lambda_u \, du \,\Big|\, \mathcal{F}_t \right), \qquad m_t^{2i} = \mathbb{E}_{\mathbb{Q}^*} \left( \int_0^{T_i} B_u^{-1} G_u \, du \,\Big|\, \mathcal{F}_t \right).$$

By the predictable representation property of the Brownian motion, $dm_t^{ji} = \psi_t^{ji} dW_t$ for $j = 1, 2$ and $i = 1, \ldots, k$ and some $\mathbb{F}$-predictable, $\mathbb{R}^d$-valued processes $\psi^{ji} = (\psi^{ji1}, \ldots, \psi^{jid})$. The proof of the next result relies on a rather straightforward verification of (45).

PROPOSITION 2.5.  *Assume that the number of traded CDS is $k = d + 1$ and the model is driven by a $d$-dimensional Brownian motion. Then conditions (45) can be represented by the linear equation $N_t \widehat{\varphi}_t = \widehat{\xi}_t$ with the $\mathbb{R}^k$-valued process $\widehat{\varphi}_t = (\varphi_t^1, \ldots, \varphi_t^k)^t$, the $\mathbb{R}^k$-valued process $\widehat{\xi}_t = (Z_t - \widetilde{S}_t, \xi_t^1, \ldots, \xi_t^d)^t$ and the $k \times k$ random matrix $N_t$ is given by*

$$N_t = \begin{bmatrix} \delta_1 \widehat{P}_t^1 - \kappa_1 \widetilde{A}_t^1 & \ldots & \delta_k \widehat{P}_t^k - \kappa_k \widetilde{A}_t^k \\ \delta_1 \psi_t^{111} - \kappa_1 \psi_t^{211} & \ldots & \delta_k \psi_t^{1k1} - \kappa_k \psi_t^{2k1} \\ \vdots & \ddots & \vdots \\ \delta_1 \psi_t^{11d} - \kappa_1 \psi_t^{21d} & \ldots & \delta_k \psi_t^{1kd} - \kappa_k \psi_t^{2kd} \end{bmatrix},$$

*where $\widehat{P}_t^i = 1 - \widetilde{P}_t^i$. For any defaultable claim $(X, A, Z, \tau)$ there exists a unique solution $(\varphi^1, \ldots, \varphi^k)$ to (45) if and only if $\det N_t \neq 0$ for almost all $t \in [0, T]$.*

**3. Multi-name credit default swap market.**  In this section we shall deal with a market model driven by a Brownian filtration in which a finite family of CDSs with different underlying names is traded.

3.1. *Price dynamics in a multi-name model.*  Our first goal is to extend the pricing results of Section 2.1 to the case of a multi-name credit risk model with stochastic default intensities.

3.1.1. *Joint survival process.*  We assume that we are given $n$ strictly positive random times $\tau_1, \ldots, \tau_n$, defined on a common probability space $(\Omega, \mathcal{G}, \mathbb{Q})$, and referred to as default times of $n$ credit names. We postulate that this space is endowed with a *reference filtration* $\mathbb{F}$, which satisfies Assumption 2.2.

In order to describe dynamic joint behavior of default times, we introduce the *conditional joint survival process* $G(u_1, \ldots, u_n; t)$ by setting, for every $u_1, \ldots, u_n, t \in \mathbb{R}_+$,

$$G(u_1, \ldots, u_n; t) = \mathbb{Q}^*(\tau_1 > u_1, \ldots, \tau_n > u_n \mid \mathcal{F}_t).$$

Let us set $\tau_{(1)} = \tau_1 \wedge \cdots \wedge \tau_n$ and let us define the process $G_{(1)}(t; t), t \in \mathbb{R}_+$ by setting

$$G_{(1)}(t; t) = G(t, \ldots, t; t) = \mathbb{Q}^*(\tau_1 > t, \ldots, \tau_n > t \mid \mathcal{F}_t) = \mathbb{Q}^*(\tau_{(1)} > t \mid \mathcal{F}_t).$$

It is easy to check that $G_{(1)}$ is a bounded supermartingale. It thus admits the unique Doob–Meyer decomposition $G_{(1)} = \mu - \nu$. We shall work throughout under the following extension of Assumption 2.1.

ASSUMPTION 3.1. We assume that the process $G_{(1)}$ is continuous and the increasing process $\nu$ is absolutely continuous with respect to the Lebesgue measure, so that $d\nu_t = v_t\,dt$ for some $\mathbb{F}$-progressively measurable, nonnegative process $v$. We denote by $\widetilde{\lambda}$ the $\mathbb{F}$-progressively measurable process defined as $\widetilde{\lambda}_t G_{(1)}^{-1}(t;t) v_t$. The process $\widetilde{\lambda}$ is called the *first-to-default intensity*.

We denote $H_t^i = \mathbb{1}_{\{\tau_i \le t\}}$ and we introduce the following filtrations $\mathbb{H}^i, \mathbb{H}$ and $\mathbb{G}$:

$$\mathcal{H}_t^i = \sigma(H_s^i; s \in [0,t]), \qquad \mathcal{H}_t = \mathcal{H}_t^1 \vee \cdots \vee \mathcal{H}_t^n, \qquad \mathcal{G}_t = \mathcal{F}_t \vee \mathcal{H}_t.$$

We assume that the usual conditions of completeness and right-continuity are satisfied by these filtrations. Arguing as in Section 2.1, we see that the process

$$\widehat{M}_t = H_t^{(1)} - \widetilde{\Lambda}_{t \wedge \tau_{(1)}} = H_t^{(1)} - \int_0^{t \wedge \tau_{(1)}} \widetilde{\lambda}_u\,du = H_t^{(1)} - \int_0^t (1 - H_u^{(1)}) \widetilde{\lambda}_u\,du$$

is a $\mathbb{G}$-martingale, where we denote $H_t^{(1)} = \mathbb{1}_{\{\tau_{(1)} \le t\}}$ and $\widetilde{\Lambda}_t = \int_0^t \widetilde{\lambda}_u\,du$. Note that the first-to-default intensity $\widetilde{\lambda}$ satisfies

$$\widetilde{\lambda}_t = \lim_{h \downarrow 0} \frac{1}{h} \frac{\mathbb{Q}^*(t < \tau_{(1)} \le t+h \mid \mathcal{F}_t)}{\mathbb{Q}^*(\tau_{(1)} > t \mid \mathcal{F}_t)} = \frac{1}{G_{(1)}(t;t)} \lim_{h \downarrow 0} \frac{1}{h} (\nu_{t+h} - \nu_t).$$

We make an additional assumption, in which we introduce the *first-to-default intensity* $\widetilde{\lambda}^i$ and the associated martingale $\widehat{M}^i$ for each credit name $i = 1, \ldots, n$.

ASSUMPTION 3.2. For any $i = 1, \ldots, n$, the process $\widetilde{\lambda}^i$ given by

$$\widetilde{\lambda}_t^i = \lim_{h \downarrow 0} \frac{1}{h} \frac{\mathbb{Q}^*(t < \tau_i \le t+h, \tau_{(1)} > t \mid \mathcal{F}_t)}{\mathbb{Q}^*(\tau_{(1)} > t \mid \mathcal{F}_t)}$$

is well defined and the process $\widehat{M}^i$, given by the formula

$$(56) \qquad \widehat{M}_t^i = H_{t \wedge \tau_{(1)}}^i - \int_0^{t \wedge \tau_{(1)}} \widetilde{\lambda}_u^i\,du,$$

is a $\mathbb{G}$-martingale.

It is worth noting that the equalities $\sum_{i=1}^n \widetilde{\lambda}^i = \widetilde{\lambda}$ and $\widehat{M} = \sum_{i=1}^n \widehat{M}^i$ are valid.

3.1.2. *Special case.* Let $\widehat{\Gamma}^i, i = 1, \ldots, n$, be a given family of $\mathbb{F}$-adapted, increasing, continuous processes, defined on a filtered probability space $(\widetilde{\Omega}, \mathbb{F}, \mathbb{P}^*)$. We postulate that $\widehat{\Gamma}^i_0 = 0$ and $\lim_{t \to \infty} \widehat{\Gamma}^i_t = \infty$. For the construction of default times satisfying Assumptions 3.1 and 3.2, we postulate that $(\widehat{\Omega}, \widehat{\mathcal{F}}, \widehat{\mathbb{P}})$ is an auxiliary probability space endowed with a family $\xi_i, i = 1, \ldots, n$, of random variables uniformly distributed on $[0,1]$ and such that their joint probability distribution is given by an $n$-dimensional copula function $C$. We then define, for every $i = 1, \ldots, n$,

$$\tau_i(\widetilde{\omega}, \widehat{\omega}) = \inf\{t \in \mathbb{R}_+ : \widehat{\Gamma}^i_t(\widetilde{\omega}) \geq -\ln \xi_i(\widehat{\omega})\}.$$

We endow the space $(\Omega, \mathcal{G}, \mathbb{Q})$ with the filtration $\mathbb{G} = \mathbb{F} \vee \mathbb{H}^1 \vee \cdots \vee \mathbb{H}^n$, where the filtration $\mathbb{H}^i$ is generated by the process $H^i_t = \mathbb{1}_{\{t \geq \tau_i\}}$ for every $i = 1, \ldots, n$.

We have that, for any $T > 0$ and arbitrary $t_1, \ldots, t_n \leq T$,

$$\mathbb{Q}^*(\tau_1 > t_1, \ldots, \tau_n > t_n \mid \mathcal{F}_T) = C(K^1_{t_1}, \ldots, K^n_{t_n}),$$

where we denote $K^i_t = e^{-\widehat{\Gamma}^i_t}$.

Schönbucher and Schubert [25] show that the following equality holds, for arbitrary $s \leq t$:

$$\mathbb{Q}^*(\tau_i > t \mid \mathcal{G}_s) = \mathbb{1}_{\{s < \tau_{(1)}\}} \mathbb{E}_{\mathbb{Q}^*} \left( \frac{C(K^1_s, \ldots, K^i_t, \ldots, K^n_s)}{C(K^1_s, \ldots, K^n_s)} \,\Big|\, \mathcal{F}_s \right).$$

Consequently, assuming that $\widehat{\Gamma}^i_t = \int_0^t \widehat{\gamma}^i_u \, du$, the $i$th survival intensity equals, on the event $\{\tau_{(1)} > t\}$,

$$\widetilde{\lambda}^i_t = \widehat{\gamma}^i_t K^i_t \frac{\frac{\partial}{\partial v_i} C(K^1_t, \ldots, K^n_t)}{C(K^1_t, \ldots, K^n_t)} = \widehat{\gamma}^i_t K^i_t \frac{\partial}{\partial v_i} \ln C(K^1_t, \ldots, K^n_t).$$

One can now easily show that the process $\widehat{M}^i$, which is given by formula (56), is a $\mathbb{G}$-martingale. This indeed follows from the Aven's lemma [1].

3.1.3. *Price dynamics of a first-to-default claim.* We will now analyze the risk-neutral valuation of first-to-default claims on a basket of $n$ credit names. As before, $\tau_1, \ldots, \tau_n$ are respective default times and $\tau_{(1)} = \tau_1 \wedge \cdots \wedge \tau_n$ stands for the moment of the first default.

DEFINITION 3.1. A *first-to-default claim* with maturity $T$ associated with $\tau_1, \ldots, \tau_n$ is a defaultable claim $(X, A, Z, \tau_{(1)})$, where $X$ is an $\mathcal{F}_T$-measurable amount payable at maturity $T$ if no default occurs prior to or at $T$, an $\mathbb{F}$-adapted, continuous process of finite variation $A : [0,T] \to \mathbb{R}$ with $A_0 = 0$ represents the dividend stream up to $\tau_{(1)}$, and $Z = (Z^1, \ldots, Z^n)$ is the vector of $\mathbb{F}$-predictable, real-valued processes, where $Z^i_{\tau_{(1)}}$ specifies the recovery received at time $\tau_{(1)}$ if default occurs prior to or at $T$ and the $i$th name is the first defaulted name, that is, on the event $\{\tau_i = \tau_{(1)} \leq T\}$.

The next definition extends Definition 2.2 to the case of a first-to-default claim. Recall that we denote $H_t^{(1)} = \mathbb{1}_{\{\tau_{(1)} \le t\}}$ for every $t \in [0, T]$.

DEFINITION 3.2. The *dividend process* $D = (D_t)_{t \in \mathbb{R}_+}$ of a first-to-default claim maturing at $T$ equals, for every $t \in \mathbb{R}_+$,

$$D_t = X \mathbb{1}_{\{T < \tau_{(1)}\}} \mathbb{1}_{[T, \infty[}(t) + \int_{]0, t \wedge T]} (1 - H_u^{(1)}) \, dA_u$$

$$+ \int_{]0, t \wedge T]} \sum_{i=1}^{n} \mathbb{1}_{\{\tau_{(1)} = \tau_i\}} Z_u^i \, dH_u^{(1)}.$$

Let us examine the price processes for the first-to-default claim. Note that

$$\mathbb{1}_{\{t < \tau_{(1)}\}} S_t^c = \mathbb{1}_{\{t < \tau_{(1)}\}} \widetilde{S}_t^c, \qquad \mathbb{1}_{\{t < \tau_{(1)}\}} S_t = \mathbb{1}_{\{t < \tau_{(1)}\}} \widetilde{S}_t,$$

where $\widetilde{S}^c$ and $\widetilde{S}$ are pre-default values of $S^c$ and $S$, where, in turn, the price processes $S^c$ and $S$ are given by Definitions 2.3 and 2.4, respectively. We postulate that, for $i = 1, \dots, n$,

$$\mathbb{E}_{\mathbb{Q}^*} |B_T^{-1} X| < \infty,$$

$$\mathbb{E}_{\mathbb{Q}^*} \left| \int_{]0, T]} B_u^{-1} (1 - H_u^{(1)}) \, dA_u \right| < \infty, \qquad \mathbb{E}_{\mathbb{Q}^*} |B_{\tau_{(1)} \wedge T}^{-1} Z_{\tau_{(1)} \wedge T}^i| < \infty,$$

so that the ex-dividend price $S_t$ (and thus also cumulative price $S^c$) is well defined for any $t \in [0, T]$. In the next auxiliary result, we denote $Y^i = B^{-1} Z^i$. Hence, $Y^i$ is a real-valued, $\mathbb{F}$-predictable process such that $\mathbb{E}_{\mathbb{Q}^*} |Y_{\tau_{(1)} \wedge T}^i| < \infty$.

LEMMA 3.1. *We have that*

$$B_t \mathbb{E}_{\mathbb{Q}^*} \left( \sum_{i=1}^{n} \mathbb{1}_{\{t < \tau_{(1)} = \tau_i \le T\}} Y_{\tau_{(1)}}^i \,\Big|\, \mathcal{G}_t \right)$$

$$= \mathbb{1}_{\{t < \tau_{(1)}\}} \frac{B_t}{G_{(1)}(t; t)} \mathbb{E}_{\mathbb{Q}^*} \left( \int_t^T \sum_{i=1}^{n} Y_u^i \widetilde{\lambda}_u^i G_{(1)}(u; u) \, du \,\Big|\, \mathcal{F}_t \right).$$

PROOF. Let us fix $i$ and let us consider the process $Y_u^i = \mathbb{1}_A \mathbb{1}_{]s, v]}(u)$ for some fixed date $t \le s < v \le T$ and some event $A \in \mathcal{F}_s$. We note that

$$\mathbb{1}_{\{s < \tau_{(1)} = \tau_i \le v\}} = H_{v \wedge \tau_{(1)}}^i - H_{s \wedge \tau_{(1)}}^i = \widehat{M}_v^i - \widehat{M}_s^i + \int_{s \wedge \tau_{(1)}}^{v \wedge \tau_{(1)}} \widetilde{\lambda}_u^i \, du.$$

Using Lemma 3.2, we thus obtain

$$\mathbb{E}_{\mathbb{Q}^*} (\mathbb{1}_{\{t < \tau_{(1)} = \tau_i \le T\}} Y_{\tau_{(1)}}^i \,|\, \mathcal{G}_t)$$

$$= \mathbb{E}_{\mathbb{Q}^*}(\mathbb{1}_A \mathbb{1}_{\{s < \tau_{(1)} = \tau_i \le v\}} \mid \mathcal{G}_t)$$

$$= \mathbb{E}_{\mathbb{Q}^*}\left( \mathbb{1}_A \left( \widehat{M}_v^i - \widehat{M}_s^i + \int_{s \wedge \tau_{(1)}}^{v \wedge \tau_{(1)}} \widetilde{\lambda}_u^i \, du \right) \Big| \mathcal{G}_t \right)$$

$$= \mathbb{E}_{\mathbb{Q}^*}\left( \mathbb{1}_A \mathbb{E}_{\mathbb{Q}^*}\left( \widehat{M}_v^i - \widehat{M}_s^i + \int_{s \wedge \tau_{(1)}}^{v \wedge \tau_{(1)}} \widetilde{\lambda}_u^i \, du \Big| \mathcal{G}_s \right) \Big| \mathcal{G}_t \right)$$

$$= \mathbb{E}_{\mathbb{Q}^*}\left( \int_{t \wedge \tau_{(1)}}^{T \wedge \tau_{(1)}} Y_u^i \widetilde{\lambda}_u^i \, du \Big| \mathcal{G}_t \right)$$

$$= \mathbb{1}_{\{t < \tau_{(1)}\}} \frac{1}{G_{(1)}(t;t)} \mathbb{E}_{\mathbb{Q}^*}\left( \int_t^T Y_u^i \widetilde{\lambda}_u^i G_{(1)}(u;u) \, du \Big| \mathcal{F}_t \right),$$

where the last equality follows from the formula

$$\mathbb{E}_{\mathbb{Q}^*}\left( \int_{t \wedge \tau_{(1)}}^{T \wedge \tau_{(1)}} R_u \, du \Big| \mathcal{G}_t \right) = \mathbb{1}_{\{t < \tau_{(1)}\}} \frac{1}{G_{(1)}(t;t)} \mathbb{E}_{\mathbb{Q}^*}\left( \int_t^T R_u G_{(1)}(u;u) \, du \Big| \mathcal{F}_t \right),$$

which is known to hold for any $\mathbb{F}$-predictable process $R$ such that the right-hand side is well defined (see Proposition 5.1.2 in [4]).  $\square$

Given Lemma 3.1, the proof of the next result is very much similar to that of Proposition 2.1 and thus is omitted.

PROPOSITION 3.1. *The pre-default ex-dividend price $\widetilde{S}$ of a first-to-default claim $(X, A, Z, \tau_{(1)})$ satisfies*

$$\widetilde{S}_t = \frac{B_t}{G_{(1)}(t;t)} \mathbb{E}_{\mathbb{Q}^*}\left( B_T^{-1} G_{(1)}(T;T) X \mathbb{1}_{\{t < T\}} \right.$$
$$\left. + \int_t^T B_u^{-1} G_{(1)}(u;u) \left( \sum_{i=1}^n Z_u^i \widetilde{\lambda}_u^i \, du + dA_u \right) \Big| \mathcal{F}_t \right).$$

By proceeding as in the proof of Proposition 2.2, one can also establish the following result, which gives dynamics of price processes $\widetilde{S}$ and $S^c$. Recall that $\mu$ is the continuous martingale arising in the Doob–Meyer decomposition of the hazard process $G_{(1)}$.

PROPOSITION 3.2. *The dynamics of the pre-default ex-dividend price $\widetilde{S}$ of a first-to-default claim $(X, A, Z, \tau_{(1)})$ on $[0, \tau_{(1)} \wedge T]$ are*

$$d\widetilde{S}_t = (r_t + \widetilde{\lambda}_t)\widetilde{S}_t \, dt - \sum_{i=1}^n \widetilde{\lambda}_t^i Z_t^i \, dt - dA_t + G_{(1)}^{-1}(t;t)(B_t \, dm_t - \widetilde{S}_t \, d\mu_t)$$
$$+ G_{(1)}^{-2}(t;t)(\widetilde{S}_t \, d\langle \mu \rangle_t - B_t \, d\langle \mu, m \rangle_t),$$

*where the continuous $\mathbb{F}$-martingale $m$ is given by the formula*

$$
\begin{aligned}
m_t = \mathbb{E}_{\mathbb{Q}^*} \bigg( & B_T^{-1} G_{(1)}(T;T) X \\
& + \int_0^T B_u^{-1} G_{(1)}(u;u) \bigg( \sum_{i=1}^n Z_u^i \widetilde{\lambda}_u^i \, du + dA_u \bigg) \, \big| \, \mathcal{F}_t \bigg).
\end{aligned}
$$
(57)

*The dynamics of the cumulative price $S^c$ on $[0, \tau_{(1)} \wedge T]$ are*

$$
\begin{aligned}
dS_t^c = \sum_{i=1}^n (Z_t^i - \widetilde{S}_{t-}) \, dM_t^i & + \bigg( r_t \widetilde{S}_t - \sum_{i=1}^n \widetilde{\lambda}_t^i Z_t^i \bigg) dt - dA_t \\
& + G_{(1)}^{-1}(t;t)(B_t \, dm_t - \widetilde{S}_t \, d\mu_t) \\
& + G_{(1)}^{-2}(t;t)(\widetilde{S}_t \, d\langle \mu \rangle_t - B_t \, d\langle \mu, m \rangle_t).
\end{aligned}
$$

3.1.4. *Hypothesis (H).* As in the single-name case, the most explicit results can be derived under an additional assumption of the immersion property of filtrations $\mathbb{F}$ and $\mathbb{G}$.

ASSUMPTION 3.3. We assume that any $\mathbb{F}$-martingale under $\mathbb{Q}^*$ is a $\mathbb{G}$-martingale under $\mathbb{Q}^*$. This also implies that Hypothesis (H) holds between $\mathbb{F}$ and $\mathbb{G}$. In particular, any $\mathbb{F}$-martingale is also a $\mathbb{G}^i$-martingale for $i = 1, \ldots, n$, that is, Hypothesis (H) holds between $\mathbb{F}$ and $\mathbb{G}^i$ for $i = 1, \ldots, n$.

It is worth stressing that, in general, there is no reason to expect that any $\mathbb{G}^i$-martingale is necessarily a $\mathbb{G}$-martingale. We shall argue that even when the reference filtration $\mathbb{F}$ is trivial this is not the case, in general (except for some special cases, e.g., under the independence assumption).

EXAMPLE 3.1. Let us take $n = 2$ and let us denote $G_t^{1|2} = \mathbb{Q}^*(\tau_1 > t \mid \mathcal{H}_t^2)$ and $G(u, v) = \mathbb{Q}(\tau_1 > u, \tau_2 > v)$. It is then easy to prove that

$$
\begin{aligned}
dG_t^{1|2} = \bigg( \frac{\partial_2 G(t,t)}{\partial_2 G(0,t)} - \frac{G(t,t)}{G(0,t)} \bigg) dM_t^2 & \\
& + \bigg( H_t^2 \partial_1 h(t, \tau_2) + (1 - H_t^2) \frac{\partial_1 G(t,t)}{G(0,t)} \bigg) dt,
\end{aligned}
$$

where $h(t, u) = \frac{\partial_2 G(t,u)}{\partial_2 G(0,u)}$ and $M^2$ is the $\mathbb{H}^2$-martingale given by

$$
M_t^2 = H_t^2 + \int_0^{t \wedge \tau_2} \frac{\partial_2 G(0,u)}{G(0,u)} \, du.
$$

If Hypothesis (H) holds between $\mathbb{H}^2$ and $\mathbb{H}^1 \vee \mathbb{H}^2$, then the martingale part in the Doob–Meyer decomposition of $G^{1|2}$ vanishes. We thus see that Hypothesis (H) is not always valid, since clearly

$$\frac{\partial_2 G(t,t)}{\partial_2 G(0,t)} - \frac{G(t,t)}{G(0,t)}$$

does not vanish, in general. One can note that in the special case when $\tau_2 < \tau_1$, the martingale part in the above-mentioned decomposition disappears and, thus, Hypothesis (H) holds between $\mathbb{H}^2$ and $\mathbb{H}^1 \vee \mathbb{H}^2$ (this case was recently studied by Ehlers and Schönbucher [14]).

From now on, we shall work under Assumption 3.3. In that case, the dynamics of price processes obtained in Proposition 3.1 simplify, as the following result shows.

COROLLARY 3.1. *The pre-default ex-dividend price $\widetilde{S}$ of a first-to-default claim $(X, A, Z, \tau_{(1)})$ satisfies*

$$d\widetilde{S}_t = (r_t + \widetilde{\lambda}_t)\widetilde{S}_t\,dt - \sum_{i=1}^{n}\widetilde{\lambda}_t^i Z_t^i\,dt - dA_t + B_t G_{(1)}^{-1}(t;t)\,dm_t,$$

*where the continuous $\mathbb{F}$-martingale $m$ is given by (57). The cumulative price $S^c$ of a first-to-default claim $(X, A, Z, \tau_{(1)})$ is given by the expression, for $t \in [0, T \wedge \tau_{(1)}]$,*

$$(58) \qquad dS_t^c = r_t S_t^c\,dt + \sum_{i=1}^{n}(Z_t^i - \widetilde{S}_t)\,d\widehat{M}_t^i + B_t G_{(1)}^{-1}(t;t)\,dm_t.$$

*Equivalently, for $t \in [0, T \wedge \tau_{(1)}]$,*

$$(59) \qquad dS_t^c = r_t S_t^c\,dt + \sum_{i=1}^{n}(Z_t^i - \widetilde{S}_t)\,d\widehat{M}_t^i + B_t G_{(1)}^{-1}(t;t)\,d\widehat{m}_t,$$

*where $\widehat{m}$ is a $\mathbb{G}$-martingale given by $\widehat{m}_t = m_{t \wedge \tau_{(1)}}$ for every $t \in [0, T]$.*

In what follows, we assume that $\mathbb{F}$ is generated by a Brownian motion. Then there exists an $\mathbb{F}$-predictable process $\xi$ for which $dm_t = \xi_t\,dW_t$ so that formula (59) yields the following result.

COROLLARY 3.2. *The discounted cumulative price of a first-to-default claim $(X, A, Z, \tau_{(1)})$ satisfies, for $t \in [0, T \wedge \tau_{(1)}]$,*

$$dS_t^{c,*} = \sum_{i=1}^{n} B_t^{-1}(Z_t^i - \widetilde{S}_t)\,d\widehat{M}_t^i + G_{(1)}^{-1}(t;t)\xi_t\,dW_t.$$

3.1.5. *Price dynamics of a CDS.* By the $i$th CDS we mean the credit default swap written on the $i$th reference name, with the maturity date $T_i$, the constant spread $\kappa_i$ and the protection process $\delta^i$, as specified by Definition 2.5.

Let $S_{t|j}^i(\kappa_i)$ stand for the ex-dividend price at time $t$ of the $i$th CDS on the event $\tau_{(1)} = \tau_j = t$ for some $j \neq i$. This value can be represented using a suitable extension of Proposition 3.1, but we decided to omit the derivation of this pricing formula. Assuming that we have already computed $S_{t|j}^i(\kappa_i)$, the $i$th CDS can be seen, on the random interval $[0, T_i \wedge \tau_{(1)}]$, as a first-to-default claim $(X, A, Z, \tau_{(1)})$ with $X = 0, Z = (S_{t|1}^i(\kappa_i), \dots, \delta^i, \dots, S_{t|n}^i(\kappa_i))$ and $A_t = -\kappa_i t$. This observation applies also to the random interval $[0, T \wedge \tau_{(1)}]$ for any fixed $T \leq T_i$.

Let us denote by $n^i$ the following $\mathbb{F}$-martingale:

$$n_t^i = \mathbb{E}_{\mathbb{Q}^*}\left( \sum_{i=1}^n \int_0^{T_i} B_u^{-1} G_{(1)}(u; u) \left( \delta_u^i \widetilde{\lambda}_u^i + \sum_{j=1, j \neq i}^n S_{u|j}^i(\kappa_i) \widetilde{\lambda}_u^j - \kappa_i \right) du \,\Big|\, \mathcal{F}_t \right).$$

The following result can be easily deduced from Proposition 3.1.

COROLLARY 3.3. *The cumulative price of the $i$th CDS satisfies, for $t \in [0, T_i \wedge \tau_{(1)}]$,*

$$dS_t^{c,i}(\kappa_i) = r_t S_t^{c,i}(\kappa_i)\, dt + (\delta_t^i - \widetilde{S}_t^i(\kappa_i))\, d\widehat{M}_t^i$$
$$+ \sum_{j=1, j \neq i}^n (S_{t|j}^i(\kappa_i) - \widetilde{S}_t^i(\kappa_i))\, d\widehat{M}_t^j + B_t G_{(1)}^{-1}(t; t)\, dn_t^i.$$

*Consequently, the discounted cumulative price of the $i$th CDS satisfies, for $t \in [0, T_i \wedge \tau_{(1)}]$,*

$$dS_t^{c,i,*}(\kappa_i) = B_t^{-1}(\delta_t^i - \widetilde{S}_t^i(\kappa_i))\, d\widehat{M}_t^i$$
$$+ B_t^{-1} \sum_{j=1, j \neq i}^n (S_{t|j}^i(\kappa_i) - \widetilde{S}_t^i(\kappa_i))\, d\widehat{M}_t^j + G_{(1)}^{-1}(t; t) \zeta_t^i\, dW_t,$$

*where $\zeta^i$ is an $\mathbb{F}$-predictable process such that $dn_t^i = \zeta_t^i\, dW_t$.*

Note that the $\mathbb{F}$-martingale $n^i$ can be replaced by the $\mathbb{G}$-martingale $\widehat{n}_t^i = n_{t \wedge \tau_{(1)}}^i$.

3.2. *Replication of a first-to-default claim.* Our final goal is to extend Theorem 3.1 in Bielecki, Jeanblanc and Rutkowski [8] and Theorem 2.1 of Section 2 to the case of several credit names in a hazard process model in which credit spreads are driven by a multi-dimensional Brownian motion. We consider a self-financing trading strategy $\varphi = (\varphi^0, \ldots, \varphi^k)$ with $\mathbb{G}$-predictable components, as defined in Section 2.2. The 0th traded asset is thus the savings account; the remaining $k$ primary assets are single-name CDSs with different underlying credit names and/or maturities. As before, for any $l = 1, \ldots, k$, we will use the short-hand notation $S^l(\kappa_l)$ and $S^{c,l}(\kappa_l)$ to denote the ex-dividend and cumulative prices of CDSs with respective dividend processes $D(\kappa_l, \delta^l, T_l, \tilde{\tau}_l)$ given by formula (25). Note that here $\tilde{\tau}_l = \tau_j$ for some $j = 1, \ldots, n$. We will thus write $\tilde{\tau}_l = \tau_{j_l}$ in what follows.

REMARK.   Note that, typically, we will have $k = n + d$ so that the number of traded assets will be equal to $n + d + 1$.

Recall that the cumulative price of a first-to-default claim $(X, A, Z, \tau_{(1)})$ is denoted as $S^c$. We adopt the following natural definition of replication of a first-to-default claim.

DEFINITION 3.3.   We say that a self-financing strategy $\varphi = (\varphi^0, \ldots, \varphi^k)$ *replicates* a first-to-default claim $(X, A, Z, \tau_{(1)})$ if its wealth process $V(\varphi)$ satisfies the equality $V_{t \wedge \tau_{(1)}}(\varphi) = S^c_{t \wedge \tau_{(1)}}$ for any $t \in [0, T]$.

When dealing with replicating strategies in the sense of the definition above, we may and do assume, without loss of generality, that the components of the process $\varphi$ are $\mathbb{F}$-predictable processes. This is rather obvious, since prior to default any $\mathbb{G}$-predictable process is equal to the unique $\mathbb{F}$-predictable process.

The following result is a counterpart of Lemma 2.3. Its proof follows easily from Lemma 2.2 combined with Corollary 3.3, and thus it is omitted.

LEMMA 3.2.   *We have, for any $t \in [0, T \wedge \tau_{(1)}]$,*

$$dV_t^*(\varphi) = \sum_{l=1}^k \varphi_t^l (B_t^{-1}(\delta_t^l - \widetilde{S}_t^l(\kappa_l))\, d\widehat{M}_t^{j_l}$$

$$+ \sum_{j=1, j \neq j_l}^n B_t^{-1}(S_{t|j}^l(\kappa_l) - \widetilde{S}_t^l(\kappa_l))\, d\widehat{M}_t^j + G_{(1)}^{-1}(t; t)\, dn_t^l),$$

*where*

$$n_t^l = \mathbb{E}_{\mathbb{Q}^*}\left(\int_0^{T_l} B_u^{-1} G_{(1)}(u; u)\left(\delta_u^l \widetilde{\lambda}_u^{j_l} + \sum_{j=1, j \neq j_l}^n S_{u|j}^l(\kappa_l)\widetilde{\lambda}_u^j - \kappa_l\right) du \,\Big|\, \mathcal{F}_t\right).$$

We are now in the position to extend Theorem 2.1 to the case of a first-to-default claim on a basket of $n$ credit names. It is also an extension of Theorem 3.1 in [8] to the case of nontrivial reference filtration $\mathbb{F}$.

Recall that $\xi$ and $\zeta^l, l = 1, \ldots, k$ are $\mathbb{F}$-predictable, $\mathbb{R}^d$-valued processes such that $dm_t = \xi_t \, dW_t$ and $dn_t^l = \zeta_t^l \, dW_t$.

THEOREM 3.1. *Assume that the processes* $\widetilde{\varphi}^1, \ldots, \widetilde{\varphi}^n$ *satisfy, for* $t \in [0, T]$ *and* $i = 1, \ldots, n$,

$$\sum_{l=1, j_l = i}^{k} \widetilde{\varphi}_t^l (\delta_t^l - \widetilde{S}_t^l(\kappa_l)) + \sum_{l=1, j_l \neq i}^{k} \widetilde{\varphi}_t^l (S_{t|i}^l(\kappa_l) - \widetilde{S}_t^l(\kappa_l)) = Z_t^i - \widetilde{S}_t$$

*and* $\sum_{l=1}^{k} \widetilde{\varphi}_t^l \zeta_t^l = \xi_t$. *Let us set* $\varphi_t^i = \widetilde{\varphi}^i(t \wedge \tau_{(1)})$ *for* $i = 1, \ldots, k$ *and* $t \in [0, T]$. *Let the process* $V(\varphi)$ *be given by Lemma 3.2 with the initial condition* $V_0(\varphi) = S_0^c$ *and let* $\varphi^0$ *be given by*

$$V_t(\varphi) = \varphi_t^0 B_t + \sum_{l=1}^{k} \varphi_t^l S_t^l(\kappa_l).$$

*Then the self-financing strategy* $\varphi = (\varphi^0, \ldots, \varphi^k)$ *replicates the first-to-default claim* $(X, A, Z, \tau_{(1)})$.

PROOF. The proof goes along the similar lines as the proof of Theorem 2.1. It suffices to examine replicating strategy on the random interval $[0, T \wedge \tau_{(1)}]$. In view of Lemma 3.2, the wealth process of a self-financing strategy $\varphi$ satisfies on $[0, T \wedge \tau_{(1)}]$

$$dV_t^*(\varphi) = \sum_{l=1}^{k} \widetilde{\varphi}_t^l (B_t^{-1}(\delta_t^l - \widetilde{S}_t^l(\kappa_l)) \, d\widehat{M}_t^{j_l}$$

$$+ \sum_{j=1, j \neq j_l}^{n} B_t^{-1}(S_{t|j}^l(\kappa_l) - \widetilde{S}_t^l(\kappa_l)) \, d\widehat{M}_t^j + G_{(1)}^{-1}(t; t)\zeta_t^l \, dW_t),$$

whereas the discounted cumulative price of a first-to-default claim $(X, A, Z, \tau_{(1)})$ satisfies on the interval $[0, T \wedge \tau_{(1)}]$ [cf. (58)]

$$dS_t^{c*} = \sum_{i=1}^{n} B_t^{-1}(Z_t^i - S_{t-}) \, d\widehat{M}_t^i + (1 - H_t^{(1)})G_{(1)}^{-1}(t; t)\xi_t \, dW_t.$$

A comparison of the last two formulae leads directly to the stated conditions. It thus suffices to verify that the strategy $\varphi = (\varphi^0, \ldots, \varphi^k)$ introduced in the statement of the theorem replicates a first-to-default claim in the sense of Definition 3.3. Since this verification is rather standard, it is left to the reader. □

T. R. Bielecki
Department of Applied Mathematics
Illinois Institute of Technology
Chicago, Illinois 60616
USA
E-mail: bielecki@iit.edu

M. Jeanblanc
Département de Mathématiques
Université d'Évry Val d'Essonne
91025 Évry Cedex
France
and
Institut Europlace de Finance
39-41 rue Cambon
F-75001 Paris
France
E-mail: mjeanbl@univ-evry.fr

M. Rutkowski
School of Mathematics and Statistics
University of New South Wales
Sydney, NSW 2052
Australia
and
Faculty of Mathematics and Information Science
Warsaw University of Technology
00-661 Warszawa
Poland
E-mail: marekr@maths.unsw.edu.au